\documentclass[11pt,leqno]{amsart}
\usepackage{amsfonts,amssymb}

\newtheorem{theorem}{Theorem}[section]
\newtheorem{lemma}[theorem]{Lemma}
\newtheorem{proposition}[theorem]{Proposition}

\newtheorem{definition}[theorem]{Definition}

\newcommand{\R}{{\mathbb R}}
\newcommand{\Z}{{\mathbb Z}}
\newcommand{\N}{{\mathbb N}}

\newcommand{\supp}{{\rm supp\ }}
\newcommand{\dd}{{\rm\; d}}
\newcommand{\ddc}{{\rm d}}

\newcommand{\trace}{{\mathrm {Tr \,}}}

\newcommand{\cL}{{\mathcal L}}

\newcommand{\LL}{\mathcal{L}}

\newcommand{\NN}{\mathcal{N}}

\newcommand{\TT}{\mathcal{T}}
\newcommand{\ND}{\rho}
\newcommand{\cP}{\mathcal{P}}

\newcommand{\cB}{\mathcal{B}}
\newcommand{\cT}{\mathcal{T}}
\newcommand{\cG}{\mathcal{G}}
\newcommand{\cK}{\mathcal{L}}

\renewcommand {\=}[1]{\stackrel{\text{#1}}{=}}

\newcommand{\cA}{{\mathcal A}}
\newcommand{\D}{{\mathrm D}}

\numberwithin{equation}{section}
\begin{document}

\title[Nonautonomous Ornstein-Uhlenbeck Equations]{Invariant Measures 
and Maximal $L^2$ Regularity for
Nonautonomous Ornstein-Uhlenbeck Equations}

\author{Matthias Geissert, Alessandra Lunardi}

\address{FB Mathematik\\ Schlossgartenstr. 7\\ TU Darmstadt\\ 64289
Darmstadt\\ Germany}
\email{geissert@mathematik.tu-darmstadt.de}

\address{Dipartimento di Matematica\\ Parco Area delle Scienze 53/A \\ 
 43100 Parma\\ Italia}
\email{alessandra.lunardi@unipr.it}


\thanks{The first author was partially supported by the MIUR - PRIN 2004 Research 
project ``Kolmogorov Equations''.}

\begin{abstract}
    We characterize the domain of the 
    realizations of the linear 
    parabolic operator $\cG$ defined by \eqref{timeoperator} in 
$L^{2}$ 
    spaces with respect to a suitable measure, that is invariant for 
    the associated evolution semigroup. As a byproduct, we obtain 
    optimal $L^{2}$ regularity results for evolution equations with 
time-depending Ornstein-Uhlenbeck operators.
\end{abstract}

\maketitle

\section{Introduction}

Finite dimensional Ornstein-Uhlenbeck operators are elliptic 
(possibly degenerate) differential operators of the type 
$$\cK \varphi (x)  = \frac12 \trace \left( B B^* \D ^2\varphi
    (x)\right)+ \langle A x ,\D \varphi(x)\rangle ,\quad x\in\R^n,$$
where $A$, $B$ are given nonzero matrices. They are prototypes of 
Kolmogorov operators,  and have been the object of several studies  
from the probabilistic and the deterministic point of view. 
The theory of linear elliptic and parabolic equations involving an 
Ornstein-Uhlenbeck operator, such as 
\begin{equation}
\lambda \varphi(x) -  \cK \varphi (x)= f(x), \quad x\in \R^n, 
\end{equation}
\begin{equation}
\label{cp}
\left\{\begin{array}{l}
u_t(t,x) + \cK u(t ,x) = g(t,x), \;\;  t>0, \;x\in \R^{n},
\\
\\
u(0,x)=\varphi(x), \;x\in \R^{n}, 
\end{array}\right.
\end{equation}
is now well developed, both in spaces of continuous bounded 
functions, and in $L^p$ spaces. See e.g. 
\cite{DPZ92,Lun97,CG01,DPZ02,MPP02,MPRS02,FL06}. The most natural 
$L^p$ setting for Ornstein-Uhlenbeck operators (as, more generally, 
for Kolmogorov operators) are not the usual $L^p$ spaces with respect 
to the Lebesgue measure, but $L^p$ spaces with respect to an 
invariant measure $\zeta$, that is a Borel probability measure in 
$\R^n$ satisfying
\begin{equation}
\int_{\R^n} T(t)\varphi \dd \zeta = \int_{\R^n} \varphi \dd \zeta , 
\quad t>0, \; \varphi \in C_b(\R^n),
 \label{invariantmeasure}
\end{equation}
where $T(t)$ is the Ornstein-Uhlenbeck semigroup.  Indeed, invariant 
measures arise naturally in the study of the asymptotic behavior of 
$T(t)$, and the realizations of $\cK$   in $L^{p}$ spaces with 
respect to invariant measures are dissipative, and therefore they 
enjoy nice analytic properties.

The aim of this paper is to extend a part of the theory to 
nonautonomous problems with time depending Ornstein-Uhlenbeck 
operators. Precisely,   we consider the  parabolic operator in 
$H^{1,2}_\mathrm{loc}(\R^{1+n})$
\begin{equation}
    \cG u(t,x)=\partial_tu(t,x)+\cK(t)u(t,\cdot)(x),\quad
    t\in\R,\ x\in\R^n ,
    \label{timeoperator}
\end{equation}
where   $\cK(t)$ are given (on $H^{ 2}_\mathrm{loc}(\R^{ n})$) by
\begin{equation}
    \cK(t)\varphi(x)=\frac12 \trace \left( B(t)B^*(t)\D_x^2\varphi
    (x)\right)+ \langle A(t)x+f(t),\D_x\varphi(x)\rangle ,\quad 
x\in\R^n,
    \label{ouoperator}
\end{equation}
with continuous and bounded data $A,B:\R\to\cL(\R^n)$ and 
$f:\R\to\R^n$.

As in the autonomous case, the operator $\cG $ arises from linear  
stochastic 
Cauchy problems in $\R^{n}$, 
\begin{equation}
\label{e1}
\left\{\begin{array}{l}
\ddc X_{t}=(A(t)X_{t}+f(t))\ddc t+B(t)\ddc W(t), \quad t\geq s, 
\\
\\
X_{s}=x,
\end{array}\right.
\end{equation}
where  $W(t)$ is a standard $n$-dimensional Brownian motion and $s\in 
\R$, $x\in \R^n$. Indeed, it is well known that, denoting by 
$X(s,t,x)$ the solution to \eqref{e1}, for each $t\in \R$ and 
$\varphi \in 
C^{2}_{b}(\R^{n})$ the function $u(s,x) :={\mathbb 
E}(\varphi(X(s,t,x)))$ satisfies the backward Kolmogorov Cauchy 
problem 
\begin{equation}
\label{e6}
\left\{\begin{array}{l}
u_s(s,x) + \cK(s)u(s,x)=0, \;\;s\leq t, \;x\in \R^{n},
\\
\\
u(t)=\varphi(x), \;x\in \R^{n}. 
\end{array}\right.
\end{equation}
See e.g. \cite{GS72,KS91}. However, our approach is purely deterministic 
and it relies on the study of the backward evolution operator 
$P_{s,t}$ for \eqref{e6} and of the associated evolution 
semigroup.

Throughout the paper  we  assume that $\cG $  is uniformly parabolic, 
i.e. there exists $\mu_0>0$ such that
\begin{equation}
    \|B(t)x\|\geq\mu_0\|x\|,\quad t\in\R,\ x\in\R^n.
    \label{invertb}
\end{equation}
Moreover, we assume that the evolution family  $U(t,s)$  generated by 
$A(\cdot)$ is stable, and we denote 
by $\omega_0(U)$ its growth bound. In other words,
\begin{equation}
\label{exp_decay}
\begin{array}{ll}
\omega_0(U) := \inf\{ & \omega \in \R :\, \exists M=M(\omega) \,{\rm  
such}\,{\rm that}
\\
\\
 & \|U(t,s)\|\leq M{\rm e}^{ \omega(t-s)},\quad -\infty<s\leq 
t<\infty \} <0.
\end{array}
\end{equation}
While in the autonomous case  these assumptions  imply existence and 
uniqueness of a probability measure $\zeta$ satisfying 
\eqref{invariantmeasure},  
in our nonautonomous case there does not exist a unique $\zeta$ such 
that 
$$\int_{\R^{n}} P_{s,t}\varphi \dd \zeta =  \int_{\R^{n}}  \varphi 
\dd 
\zeta, \quad s<t, \;\varphi \in C_{b}(\R^{n}),$$
but we can find families  of measures $\{\mu_{t}:$ $t\in \R\}$, 
called {\em entrance laws at time} $-\infty$ in \cite{Dyn89} and 
{\em evolution systems of measures} in \cite{DPR05},  such that 
\begin{equation}
\label{evmeas}
\int_{ \R^n} P_{s,t}\varphi  \dd   \mu_s = 
\int_{ \R^n} \varphi  \dd   \mu_t, \quad \varphi \in 
C_{b}(\R^{n}), \;s\leq t.
\end{equation}
Such families of measures are infinitely many, and they are 
characterized in Section 2. 
Among all of them, the simplest one consists of the Gaussian measures 
$\nu_t$ defined by
\begin{equation}
\label{nu_t}
 \nu_t=\NN_{g(t,-\infty),Q(t,-\infty)},\quad t\in\R,
\end{equation}
where%
$$g(t,s):=\int\limits_s^tU(t,r)f(r)\dd r,\quad
Q(t,s):=\int\limits_s^tU(t,r)B(r)B^*(r)U^*(t,r)\dd r, \quad -\infty 
\leq s < t, $$
and we denote by  $\NN_{m,Q}$ the $n$-dimensional Gaussian
measure with covariance operator $Q$ and mean $m$.

With the aid of the measures $ \nu_t$,   we construct  a Borel 
measure in $ \R^{1+n}$ as follows: for $I\in\cB(\R)$ and 
$K\in\cB(\R^n)$, we set  $\nu(I\times K)=\int_I\nu_t(K)\dd t$, then 
$\nu$ is extended in a standard way to a measure on $\cB(\R^{1+n})$, 
still denoted by $\nu$.
 
We define $G_0:D(G_{0})\subset L^2(\R^{1+n},\nu)\to 
L^2(\R^{1+n},\nu)$ by
\begin{align*}
    \left(G_0u\right)(t,x)&=\cG u(t,x),\quad
    x\in\R^n,\ t\in \R,\ u\in D(G_0), 
\end{align*}
where $D(G_0)$ is a core of nice functions; precisely, it is the 
linear span 
of real and imaginary parts of the functions $u$ of the type $u(t,x) 
=\Phi_j(t) e^{i\langle x,h_j(t)\rangle}$ with 
$\Phi_j\in C_c^1(\R)$ and $h_j\in C^1_b(\R;\R^n)$.  
Then $\nu$ is invariant for $G_{0}$, in the sense that
\begin{equation}
    \int_{\R\times\R^n} G_{0} u(t,x) \ddc \nu = 0, \quad u\in 
D(G_{0}).
\end{equation}
A fundamental property of the realizations of second order elliptic 
and 
parabolic operators in $L^{2}$ spaces with respect to invariant 
measures is their dissipativity. In fact, since $G_0(u^{2}) = 
2u\,G_{0}u + |B^{*} D_{x}u|^{2}$ and the integral of $G_0(u^{2})$ 
vanishes, we get
\begin{align*}
   \int_{\R^{1+n}} u\,G_{0}u \dd\nu = -\frac{1}{2} 
   \int_{\R^{1+n}}|B^{*} D_{x}u|^{2} \dd\nu\leq 0, 
\end{align*}
so that $\langle u, G_{0}u\rangle_{L^{2}}\leq 0$ for each $u\in 
D(G_{0})$. Being dissipative,  $G_0$ is closable. Its closure  
 $G$ is dissipative and has dense domain because $D(G_{0})$ is dense.  
 $G$ is the natural realization of 
$\cG$ in $L^2(\R^{1+n},\nu)$, as next theorem \ref{main1} states.

For $k$, $s\in \N$ we define the Sobolev spaces
\begin{align*}
    H^{k,s}(\R^{1+n},\nu):=\big\{u\in
    H^{k,s}_\mathrm{loc}(\R^{1+n}):&\;\partial_t^lu\in
    L^2(\R^{1+n},\nu)\mbox{ for all }0<l<k,\\& \;\D_x^\alpha u\in
    L^2(\R^{1+n},\nu)\mbox{ for all }|\alpha|\leq s\big\}.
\end{align*}
Our first main result reads as follows. 
\begin{theorem}\label{main1}
    We have 
    \begin{align*}
        D(G)=H^{1,2}(\R^{1+n},\nu)=\left\{u\in
        H^{1,2}_\mathrm{loc}(\R^{1+n})\cap
                L^2(\R^{1+n},\nu):\cG
                                u\in
         L^2(\R^{1+n},\nu)\right\}.
    \end{align*}
\end{theorem}

Note that  $\nu$   is not a probability measure, because of the 
Lebesgue measure with respect to time in the whole $\R$. 
To avoid this drawback we may  work in spaces of time periodic 
functions, assuming that also the coefficients $A$, $B$, $f$ are 
periodic with the same period $T$. Then $\frac{1}{T}\nu$ is a 
probability measure in $(0,T)\times \R^{n}$. 

To be precise, let $L^2_\#( \R^{1+n})$ denote the Hilbert 
space 
of all Lebesgue measurable functions $u:\R^{1+n}\to\R$ such that 
$u(t,x)=u(t+T,x)$
a.e. $t\in\R$, $x\in\R^n$ and $u|_{(0,T)\times\R^n}\in L^2(
(0,T)\times\R^n,\nu)$, endowed with the norm 
$$u\mapsto \bigg( \frac{1}{T} \int_{0}^{T}\int_{\R^{n}} u(t,x)^{2} 
\ddc \nu_{t}\dd t\bigg)^{1/2}.$$
Similarly as above, and as in \cite{DPL06}, we define 
$G_0^\#:D(G_0^\#)\subset L^2_\#(\R^{1+n},\nu)\mapsto 
L^2_\#(\R^{1+n},\nu)$ 
where $D(G_{0}^\#)$ is the linear span 
of real and imaginary parts of the functions $u$ of the type $u(t,x) 
=\Phi_j(t) e^{i\langle x,h_j(t)\rangle}$ with  $T$-periodic 
$\Phi_j\in C^1 (\R)$ and $h_j\in C^1 (\R;\R^n)$. 
Again,  $G_0^\#$ is dissipative, hence closable, and its closure 
$G^\#$ generates a $C_0$-semigroup $(\cP^\#_\tau)_{\tau\geq0}$ of 
contractions on 
$L^2_\#(\R^{1+n},\nu)$. 

For $k$, $s\in \N$ we set 
\begin{align*}
    H_\#^{k,s}(\R^{1+n},\nu):=\big\{u\in
    H^{k,s}_\mathrm{loc}(\R^{1+n}):&\;\partial_t^lu\in
    L^2_\#(\R^{1+n},\nu)\mbox{ for all }0<l<k,\\& \;\D_x^\alpha u\in
    L^2_\#(\R^{1+n},\nu)\mbox{ for all }|\alpha|\leq s\big\}.
\end{align*}
The description of the domain in the $T$-periodic case reads as 
follows.

\begin{theorem}\label{sharpmain1}
We have 
$$ \begin{array}{lll}
        D(G_\#)& =& H^{1,2}_\#(\R^{1+n},\nu)
	\\& =& \left\{u\in
        H^{1,2}_\mathrm{loc}(\R^{1+n})\cap
                L^2_\#(\R^{1+n},\nu):\cG
                                u\in
         L^2_\#(\R^{1+n},\nu)\right\}.
    \end{array}$$
\end{theorem}
 This characterization yields that $D(G_\#)$ is 
compactly embedded in $L^2_\#(\R^{1+n},\nu)$, \linebreak through the 
compactness of the embedding 
$H^{1,2}_\#(\R^{1+n},\nu)\subset 
L^2_\#(\R^{1+n},\nu)$  that we prove in Section 5. 

Theorems 1.1 and 1.2 may be seen as maximal regularity results for 
evolution equations with time in $\R$, 
\begin{equation}
    \label{timeinR}
    u_{t}(t,x)+\cK (t)u (t,x)= h(t,x), \quad t\in \R, \;x\in \R^{n},
\end{equation}
with datum $h\in  L^2 ( \R^{1+n},\nu)$ (respectively, $h\in  
L^2_\#(\R^{1+n},\nu)$ in the periodic case), since they 
state that if $u\in L^2 ( \R^{1+n},\nu)\cap 
H^{1,2}_\mathrm{loc}(\R^{1+n},\nu)$ (resp.,  $u\in  
L^2_\#(\R^{1+n},\nu)\cap H^{1,2}_\mathrm{loc}(\R^{1+n},\nu)$)
satisfies \eqref{timeinR} then $u_{t}$ and each second order space 
derivative 
$D_{ij}u$ belong to $ L^2 ( \R^{1+n},\nu)$ (respectively, to 
$L^2_\#(\R^{1+n},\nu)$). 

Concerning solvability of problem \eqref{timeinR}, we 
remark that it is not a Cauchy problem and we do not expect existence 
and uniqueness of a solution $u$ for any $h$; in fact, it is not hard 
to see 
that $0$ is in the spectrum of $G$ and of $G_\#$. (The spectral 
properties of $G$ and of $G_\#$, as well as asymptotic behavior of 
$P_{s,t}$, will be studied in a forthcoming 
paper \cite{GL07}). 

Note that problem \eqref{timeinR} cannot be seen as an evolution 
equation of the type $u'(t) + L(t)u(t) = h(t)$ in a fixed Hilbert 
space $H$, because the Hilbert spaces $L^2 ( \R^{ n},\nu_{t})$ 
involved here vary with time. So we cannot use the techniques of 
evolution equations in (fixed) Hilbert spaces. 

Our procedure is the following: we use the fact that $G$ and $G_\#$ 
are the infinitesimal generators 
of the evolution semigroup 
\begin{align}\label{rep}
    \left( \cP_\tau u \right)(t,x)&=\left( P_{t,t+\tau}
    u(t+\tau,\cdot) \right)(x),\quad t\in\R,\ x\in\R^n,\
    \tau\geq0,
\end{align}
in the spaces $L^{2}(\R^{1+n},\nu)$ and $L^2_\#(\R^{1+n},\nu)$, 
respectively. 
We prove optimal blow-up estimates for the space derivatives of 
$P_{s,t}\varphi$ 
near $t=s$ for any $\varphi\in L^2(\R^{n},\nu_{t})$ and any 
multi-index $\alpha$,  
\begin{equation}
    \label{smoothing}
\|\D_x^\alpha P_{s,t}\|_{\cL(L^2(\R^{n},\nu_t),L^2(\R^{n},\nu_s))}\leq
    C(t-s)^{-|\alpha|/2} ,\quad 0<t-s<1, 
\end{equation}
that yield optimal estimates for 
the norm of $\cP_\tau $ in $\cL(L^2 (\R^{1+n},\nu), 
H^{0,k}(\R^{1+n},\nu))$
and in $\cL(L^{2}_{\#} (\R^{1+n},\nu), H^{0,k}_{\#}(\R^{1+n},\nu))$ 
for 
$k\in \N$, near $\tau =0$. Then we use an interpolation 
theorem that gives optimal embeddings for the domain of the 
infinitesimal 
generator of a semigroup $T(\tau)$ from optimal estimates on the 
behavior of $T(\tau)$ near $\tau=0$; in our case it gives 
$D(L)\subset ( L^2 (\R^{1+n},\nu), 
H^{0,4}(\R^{1+n},\nu))_{1/2,2}$. The latter space is readily characterized as 
$H^{0,2}(\R^{1+n},\nu)$. 
 
The crucial step are the smoothing estimates  \eqref{smoothing} that 
are quite 
similar to the corresponding estimates in the autonomous case. 
Together with \eqref{smoothing}  we obtain also optimal estimates for 
$t-s\to \infty$, 
\begin{equation}
 \label{smoothing2}
\|\D_x^\alpha P_{s,t}\|_{\cL(L^2(\R^{n},\nu_t),L^2(\R^{n},\nu_s))}\leq
 C{\rm e}^{ \omega|\alpha|(t-s)},\quad t-s>1, 
\end{equation}
where $\omega$ is any number in $(\omega_{0}(U),0)$ 
and $C=C(\alpha, \omega)$. These estimates will be the starting point 
of the study of spectral properties and asymptotic behavior of the 
forthcoming paper \cite{GL07}.

The characterization of the domain of  $G$ in Theorem 1.1 allows us
to study maximal $L^2$ regularity in backward Cauchy problems
such as 
\begin{equation}
\label{e7}
\left\{\begin{array}{l}
u_s(s,x) + \cK(s)u(s,x)=h(s), \;\;s\in(T_1,T_2), \;x\in \R^{n},
\\
\\
u(T_2,x)=\varphi(x), \;x\in \R^{n}. 
\end{array}\right.
\end{equation}
with fixed $T_1<T_2$. Note that we do not assume the
coefficients $A$, $B$ and $f$ to be $(T_2-T_1)$-periodic.

\begin{theorem}
\label{thm:maxreg}
Let $T_1<T_2$. 
For each $h\in L^2( (T_1,T_2)\times\R^n,\nu)$ and $\varphi \in
H^1(\R^n,\nu_{T_2})$ there exists a unique solution $u\in
H^{1,2}( (T_1,T_2)\times\R^n,\nu)$ of \eqref{e7}. Moreover, $u$ 
satisfies
\begin{align}
\label{estmaxreg}
\|u\|_{H^{1,2}( (T_1,T_2)\times\R^n,\nu)}\leq C\left(
\|h\|_{L^2( (T_1,T_2)\times\R^n,\nu)}
+\|\varphi\|_{ H^1(\R^n,\nu_{T_2})}
\right),
\end{align}
where $C>0$ is independent of $h$ and $\varphi$.
\end{theorem}
The assumption $u_{0}\in H^1(\R^n,\nu_{T_2})$ is necessary for $u\in
H^{1,2}( (T_1,T_2)\times\R^n,\nu)$ because  
$H^1(\R^n,\nu_{T_2})$ is the space of the traces  $u(T_{2}, \cdot)$ 
of the 
functions   $u\in H^{1,2}( (T_1,T_2)\times\R^n,\nu)$. 

This paper is organized as follows.
In section 2 we characterize all the families of probability measures 
$\{\mu_{t}:$ $t\in \R\}$ such that 
$$\int_{ \R^n} P_{s,t}\varphi  \dd \mu_s = 
  \int_{ \R^n} \varphi  \dd \mu_t, \quad \varphi \in 
  C_{b}(\R^{n}), \;s\leq t,$$
and we show that the measures $\nu_{t}$ defined in \eqref{nu_t} are 
the
unique ones satisfying
$$ \sup_{t\in \R} \int\limits_{\R^n} |x|^{\alpha}\mu_t(\ddc x) 
<+\infty $$
for some $\alpha >0$. 
The proofs of the domain characterizations  are given in 
Sections~\ref{secmain1} and
\ref{secmain2}, respectively. Estimates \eqref{smoothing} and 
\eqref{smoothing2} are proved 
in Section 3.  The 
characterization of real interpolation spaces between 
$L^2(\R^{1+n},\nu)$
and $H^{0,k}(\R^{1+n},\nu)$ and  between
$L_\#^2(  \R^{1+n},\nu)$ and $H_\#^{0,k}(
 \R^{1+n},\nu)$ is given in 
Section~\ref{secinterpol}. 
Finally, the proof of Theorem~\ref{thm:maxreg} is given in Section 
6.


\section{Preliminaries and notation, invariant measures}


We recall some general facts about time-dependent Ornstein-Uhlenbeck 
operators and their invariant measures, already partly discussed 
in \cite{DPL06}.

First of all, for all $\varphi\in C^2_b(\R^n)$ and $t\in \R$ there 
exists a unique bounded solution $u\in C^{1,2}(\{(s,x)\in \R^{1+n}:\, 
s\leq t\})$ of the problem

\begin{equation}
    \left\{
    \begin{array}{rll}
    \partial_su(s,x)+\cK(s)u(s,x)=&0,&\quad x\in\R^n,\ s<t, \\
    u(t,x)=&\varphi(x),&\quad x\in\R^n.
    \end{array}
    \right.
    \label{ivp}
\end{equation}
The solution $u(s,x):=P_{s,t}\varphi(x)$ of \eqref{ivp} is 
given by the formula
\begin{align}
P_{s,t}\varphi(x)=\int\limits_{\R^n}\varphi\left(
y+g(t,s)
\right)\NN_{U(t,s)x,Q(t,s)}(\ddc y),\quad
-\infty<s\leq t<\infty.
\label{evolution_operator}
\end{align}

This has been shown in \cite{DPL06} under periodicity assumptions on 
the coefficients, but the proof goes through as well in the general 
case. It is easy to see 
that under our non-degeneration assumption \eqref{invertb}, for each 
$\varphi\in C_b(\R^n)$ the function  $u(s,x) 
= P_{s,t}\varphi(x)$ is still the unique bounded classical solution 
to 
\eqref{ivp}. The associated  evolution semigroup in $C_{b}(\R^{1+n})$ 
is defined by 
$$\cP_\tau u(t,x) = P_{t,t+\tau}u(t+\tau, \cdot)(x), \quad \tau\geq0,\
x\in\R^n,\ t\in\R.$$

\begin{definition}
    \label{definv}
A measure $\nu$ in $\R^{1+n}$ is said to be invariant for $ 
\cP_\tau $ if 
\begin{equation}
    \int\limits_{\R^{1+n}}\cP_\tau u  \dd \nu = 
\int\limits_{\R^{1+n}}  u 
 \dd \nu , \quad \tau >0, \;u\in C_{b}(\R^{1+n})\cap L^{1}(\R^{1+n}, 
\nu). 
\label{misurainvariante}
\end{equation}

An evolution system of measures for 
$P_{s,t}$ is a family of probability measures $(\nu_t)_{t\in\R}$ in 
$\R^{n}$ such that 
\begin{equation}
    \int\limits_{\R^n}P_{s,t}\varphi(x)\nu_s(\ddc x)=
    \int\limits_{\R^n}\varphi(x)\nu_t(\ddc x),
    \quad \varphi\in
    C_b(\R^n).
    \label{invariantmess}
\end{equation}
\end{definition}

If $P_{s,t}$ has an  evolution system of 
measures $(\nu_t)_{t\in\R}$, then one can define an invariant measure 
for $ \cP_\tau $, as follows. For Borel sets $I\subset \R$, $K\subset 
\R^{n}$ we define $\nu(I\times K)=\int_I\nu_t(K)\dd t$, then $\nu $ 
is extended in a standard 
way to all $\cB( \R^{1+n})$. 
It is easy to see that  $\mu$ is invariant for $\cP_\tau $.

It is well known (Khas'minski\u{i} Theorem) that  if a Markov 
semigroup is strong Feller and 
irreducible, then it has at most one invariant measure.  Our 
semigroup $\cP_\tau $ is not irreducible and it
is not strong Feller because of the translation part, and it has in 
fact infinitely many invariant measures, as the next proposition 
shows. 
For 
its proof, we recall that the Fourier transform of a Gaussian 
measure ${\mathcal N}_{m,Q}$ is given by 
\begin{equation}
\label{Fourier}
\widehat{\mathcal N} _{m,Q}(h) = e^{i\langle m,h  
\rangle-\frac12\;\langle  
Qh,h  
\rangle}, \quad h\in \R^n. 
\end{equation}

\begin{proposition} 
    \label{infinvmeas}
Fixed any $t_{0}\in \R$ and any Borel probability measure $\mu $  in 
$\R^{n}$, define a 
family of Borel probability measures $\mu_{t}$ by their 
Fourier transforms, 
\begin{equation}
\label{defFT}
\widehat{\mu}_{t}(h) :=  \widehat{\mu} (U^{*}(t,t_{0})h), \quad t\in 
\R .
\end{equation}
Then the measures $\nu_{t}$ defined by 
\begin{equation}
    \label{zetat}
    \nu_{t} = \NN_{g(t,-\infty),Q(t,-\infty)}\star \mu_{t}, \quad 
t\in  \R ,
\end{equation}
are an evolution system of measures for $P_{s,t}$. Moreover, 
all the evolution systems of measures for $P_{s,t}$ are of this type. 
\end{proposition}
\begin{proof}
We remark that a family of Borel probability measures $( 
\nu_{t})_{t\in \R}$ is an evolution system of measures iff their 
Fourier transforms satisfy
\begin{equation}
\label{FT}
e^{ i\langle g(t,s),h\rangle} e^{-\frac{1}{2}\langle Q(t,s)h,h\rangle 
} 
\widehat{\nu}_{s}(U^{*}(t,s)h) = \widehat{\nu}_{t}(h) , \quad s\leq 
t, \;h\in \R^{n}. 
\end{equation}
Indeed,  the left hand side of \eqref{FT} is  equal to
$ \int_{\R^n}P_{s,t}\varphi(x)\nu_s(\ddc x)$ and the right 
hand side is $ \int_{\R^n}\varphi(x)\nu_t(\ddc x)$ if we take
$\varphi(x) := e^{i\langle x,h\rangle}$. So, if $( 
\nu_{t})_{t\in \R}$ is an evolution system of measures the left and 
the right hand side have to 
coincide. Conversely, 
if $\varphi$ is any continuous bounded function, there
exists a  sequence $(\varphi_{k})_{k\in {\mathbb N}}$, where
$\varphi_k\in\mathrm{span} \{e^{i\langle x,h\rangle}:h\in\R^n\}$ such 
that  $\|\varphi_{k}\|_{ \infty }\le \|\varphi \|_{ \infty}$,
and $ \lim_{k\to \infty }\varphi_{k}(x)=\varphi(x)$, 
$\forall$ $x\in  \R^n$. 
By \eqref{FT}, the equality 
$ \int_{\R^n}P_{s,t}\varphi_{k}(x)\nu_s(\ddc x)$ $=$
$ \int_{\R^n}\varphi_{k}(x)\nu_t(\ddc x)$
holds for each $k$, 
and letting $k\to \infty$ we get \eqref{invariantmess}.

If $\nu_{t}$ is defined by \eqref{zetat}, then for each $h\in \R^{n}$ 
and $s<t$ we have
$$\begin{array}{l}
e^{i\langle  g(t,s),h \rangle-\frac12\;\langle  
Q(t,s)h,h  \rangle} \widehat{\nu}_{s}(U^{*}(t,s)h) 
\\
\\
= e^{i\langle  g(t,s),h \rangle-\frac12\;\langle  
Q(t,s)h,h  \rangle} \widehat{\mathcal N}_{g(s,-\infty), 
Q(s,-\infty)}(U^{*}(t,s)h) 
\widehat{\mu}_{s}(U^{*}(t,s)h) 
\\
\\
= e^{i\langle  g(t,s),h \rangle-\frac12\;\langle Q(t,s)h,h \rangle}
e^{i\langle  g(s,-\infty),U^{*}(t,s)h \rangle-\frac12\;\langle Q( 
s,-\infty)U^{*}(t,s)h,U^{*}(t,s)h \rangle}\widehat{\mu} ( 
U^{*}(s,t_{0})  U^{*}(t,s)h) 
\\
\\
= e^{i\langle  g(t,-\infty),h \rangle-\frac12\;\langle 
Q(t,-\infty)h,h \rangle}
\widehat{\mu} (U^{*}(t,t_{0})h) 
\\
\\
=  \widehat{\nu}_{t}(h),
\end{array}$$
so that \eqref{FT} holds.

Conversely, if  \eqref{FT} holds, then the left hand side is 
independent of $s$, hence for each $t\in \R$ and $h\in 
\R^{n}$ there exists the limit 
\begin{equation}
    \label{mu0}
\lim_{s\to-\infty} \widehat{\nu}_{s}(U^{*}(t,s)h) = 
\widehat{\nu}_{t}(h) 
e^{- i\langle g(t,-\infty),h\rangle +\frac{1}{2}\langle 
Q(t,-\infty)h,h\rangle } .
\end{equation}
Being the pointwise limit of Fourier transforms of probability 
measures, by the Bochner Theorem the left hand side of \eqref{mu0} is 
the Fourier 
transform of a probability measure $\mu_{t}$, and for each $t$, 
$t_{0}\in 
\R$, $h\in \R^{n}$ we have 
$$\widehat{\mu}_{t_{0}}(U^{*}(t,t_{0})h) = \lim_{s\to -\infty} 
\widehat{\nu}_{s}(U^{*}(t_{0},s)U^{*}(t,t_{0} )h) =  
\widehat{\mu}_{t}(h)$$
because $U^{*}(t_{0},s)U^{*}(t,t_{0})h =  U^{*}(t,s)h$.  Therefore, 
$\widehat{\mu}_{t}$ satisfies \eqref{defFT} with $\mu = \mu_{t_{0}}$. 
Now \eqref{mu0} implies 
$$\widehat{\nu}_{t}(h) = e^{i\langle g(t,-\infty),h\rangle 
-\frac{1}{2}\langle 
Q(t,-\infty)h,h\rangle } \widehat{\mu}_{t}( h) ,$$
hence $\nu_{t} = \NN_{g(t,-\infty),Q(t,-\infty)}\star \mu_{t}$. 
\end{proof}

The family of measures $({\mathcal N}_{g(t,-\infty), 
Q(t,-\infty)})_{t\in \R}$ corresponds to $t_{0}=0$ and 
$\mu_{0}=\delta_{0}$. For a similar result in a much more general setting (but with 
$f\equiv 0$) see \cite[Thm.~5.1]{Dyn89}. 
In the case of $T$-periodic coefficients, it is the unique 
$T$-periodic evolution system of measures for $P_{s,t}$, see 
\cite{DPL06}. In the general case, we have the following 
characterization.

\begin{lemma} 
    \label{uniqueness}
Let $(\nu_{t})_{t\in \R}$ be an evolution system of measures for 
$P_{s,t}$ such that 
$$\exists \alpha >0:\quad \sup_{t\in \R} \int\limits_{\R^n} 
|x|^{\alpha}\nu_t(\ddc x) <+\infty.$$
Then $\nu_{t} = \NN_{g(t,-\infty),Q(t,-\infty)} $, for 
each $ t\in \R$.   
\end{lemma}
\begin{proof}
Since, by assumption,  $\int_{\R^{n}} \varphi(x) 
\nu_t(dx) = \int_{\R^{n}} P_{s,t}\varphi(x) 
\nu_s(dx)$ for each  $s\leq t$, it is enough to show that 
$$\lim_{s\to -\infty}\int_{\R^{n}} P_{s,t}\varphi(x) 
\nu_s(dx) = \int_{\R^{n}} \varphi(x) 
\NN_{g(t,-\infty),Q(t,-\infty)} (dx), \quad t\in \R, \;\;\varphi\in 
C^{1}_{b}(\R^{n}). $$
We have
\begin{align*}
\int_{\R^{n}}P_{s,t} \varphi(x) \nu_s(\ddc x)  
 =&   \int_{\R^{n}} \varphi(y)\NN_{g(t,-\infty),Q(t,-\infty)}(\ddc 
y)\\
&+ 
\int_{\R^{n}} \bigg( 
P_{s,t}\varphi(x) -  \int_{\R^{n}} \varphi(y) 
\NN_{g(t,-\infty),Q(t,-\infty)}(dy)\bigg)\nu_s(\ddc x).
\end{align*}
To prove that the last integral goes to zero as $s\to -\infty$ we 
estimate the integrand
   \begin{equation}
    \begin{array}{l}
 \left|P_{s,t}\varphi(x) - \displaystyle { \int_{\R^{n}} \varphi(y) 
       \NN_{g(t,-\infty),Q(t,-\infty)}(\ddc y) }
\right|  
\\
\\
= \left| \displaystyle {\int_{\R^{n}} (P_{s,t }\varphi( x) - P_{s,t 
}\varphi(y))\NN_{g(s,-\infty),Q(s,-\infty)}(\ddc y)}
\right|. \end{array}
\label{I}
   \end{equation}
Without loss of generality we may assume that $\alpha<1$. Fix $\omega 
\in (\omega_{0}(U), 0)$. 
Recalling that 
$$[\psi]_{C^{\alpha}(\R^{n})} \leq (2\|\psi\|_{\infty})^{ 1-\alpha} 
\| \,|D \psi|\,\|_{\infty}^{\alpha}, \quad \psi\in 
C^{1}_{b}(\R^{n}), $$
and that 
$$ | D_{x}P_{s,t}\varphi (x)|  = |U^{*}(t,s)( 
P_{s,t}(D \varphi))(x)| \leq Me^{ \omega 
(t-s)}\|\,|D \varphi|\,\|_{\infty}  \;\; s\leq t, \;x\in \R^{n}, $$
 we get 
$$|P_{s,t }\varphi( x) - P_{s,t }\varphi(y)| \leq 
(2\|\varphi\|_{\infty})^{ 1-\alpha} 
(Me^{ \omega 
(t-s)}\|\,|D \varphi|\,\|_{\infty}|x-y|)^{\alpha } := 
C_{1}e^{ \alpha \omega 
(t-s)}|x-y|^{\alpha }, $$
with $C_{1}$ independent of $t$, $s$, $x$, $y$. Therefore, 
$$\begin{array}{l}
\left|  \displaystyle {\int_{\R^{n}} (P_{s,t }\varphi( x) - P_{s,t 
}\varphi(y))\NN_{g(s,-\infty),Q(s,-\infty)}(\ddc y)}
\right| 
\\
\\
\leq C_{1}e^{\alpha \omega 
(t-s)} \displaystyle 
{\int_{\R^{n}}|x-y|^{\alpha}\NN_{g(s,-\infty),Q(s,-\infty)}(\ddc y)}. 
\end{array}$$
By the H\"older inequality  we have
\begin{align*}
    \int_{\R^{n}}|x-y|^{\alpha}\NN_{g(s,-\infty),Q(s,-\infty)}(\ddc 
y)\leq 
\left(\int_{\R^{n}}|x-y|^{2}\NN_{g(s,-\infty),Q(s,-\infty)}(\ddc 
y)\right)^{\alpha /2}
\\
= (|x-g(s,-\infty)|^{2} + \mbox{\rm Tr}\,Q(s, -\infty))^{\alpha /2}
\leq C_{2}(|x|^{\alpha}+1), 
\end{align*}
with  $C_{2}$ independent of $s$, $x$. Replacing in \eqref{I}, we get 
$$ \left|P_{s,t}\varphi(x) -  \int_{\R^{n}} 
\varphi(y)\NN_{g(t,-\infty),Q(t,-\infty)}(\ddc y) 
\right|  
\leq 
C_{1}C_{2}e^{ \alpha \omega (t-s)} ( |x|^{\alpha} +1)$$
and since $\omega <0$  the  statement follows. \end{proof}

From now on, we shall consider only the  evolution system of 
measures 
$$\nu_{t} := \NN_{g(t,-\infty),Q(t,-\infty)}$$ 
and the 
corresponding invariant measure $\nu$ for $\cP_\tau $. Note that they 
satisfy
$$\sup_{t\in \R} \int\limits_{\R^n} 
|x|^{\alpha}\nu_t(\ddc x) <+\infty, \quad \forall \alpha >0.$$

In contrast to the autonomous case, we cannot expect
\begin{align*}
\int\limits_{\R^n}\cL(s)\varphi\,\nu_s(\ddc x)=0,\quad \forall \varphi\in
H^2(\R^n,\nu_s), s \in\R.
\end{align*}
However, we have the following Lemma.
For 
its proof, we recall 
that in the nondegenerate case $\det Q\neq 0$  the density 
$\rho_{m,Q}$ of $\NN_{m,Q}$ with respect to the Lebesgue measure is given by
\begin{align*}
\rho_{m,Q}(y) = (2\pi)^{-\frac n2}(\det Q)^{-\frac12}{\rm
  e}^{-\frac12 \langle Q^{-1}(y-m) ,y-m\rangle},\quad y\in\R^n.
\end{align*}

\begin{lemma}
\label{invop}
Let $\varphi\in H^2(\R^n,\nu_s)$ for some $s\in\R$. Then,
\begin{align*}
\int\limits_{\R^n}\cL(s)\varphi\, \nu_s(\ddc x)	
=\int\limits_{\R^n}\varphi\, \partial_s\rho(s,x)\dd x,
\end{align*}
where $\rho(s,x)=\rho_{g(s,-\infty),Q(s,-\infty)}$.
\end{lemma}
\begin{proof}
Since $\mathrm{span}\{ {\rm e}^{i\langle k,\cdot\rangle}:k\in\R^n\}$
is dense in $H^2(\R^n,\nu_s)$, it suffices to show that
\begin{align*}
\int\limits_{\R^n}\cL(s){\rm e}^{i\langle
k,x\rangle}\nu_s(\ddc x)
=\int\limits_{\R^n}{\rm e}^{i\langle k,x\rangle}
\partial_s\rho(s,x)\dd x,\quad k\in\R^n.
\end{align*}
We have
\begin{align*}
	\int\limits_{\R^n}\cL(s){\rm e}^{i\langle
	k,x\rangle}\nu_s(\ddc x)
	=&\int\limits_{\R^n}\left(-\frac12|B^*
	k|^2 +i \langle A(s)x,k\rangle+i\langle
	f(s),k\rangle\right){\rm
	e}^{i\langle k,x\rangle}\nu_s(\ddc x)
	\\
	=&-\frac 12| B^* k|^2\widehat{\mathcal N}_{g(s,-\infty),Q(s,-\infty)}(k)
	+\langle A(s)\nabla \widehat{\mathcal N}_{g(s,-\infty),Q(s,-\infty)}(k),
	k \rangle
	\\
	&+i\langle f(s),k\rangle \widehat{\mathcal N}_{g(s,-\infty),Q(s,-\infty)}(k)\\
	=&\big[ -\frac 12 |B^* k|^2 
	+\langle
	A(s)\left(ig(s,-\infty)-Q(s,-\infty)k\right) ,k\rangle\\
	&+i\langle f(s),k\rangle\big]\widehat{\mathcal N}_{g(s,-\infty),Q(s,-\infty)}(k)\\
	&=\partial_s\widehat{\mathcal N}_{g(s,-\infty),Q(s,-\infty)}(k)
	=\int\limits_{\R^n}{\rm e}^{i\langle k,x\rangle}
	\partial_s\rho(s,x)\dd x.
	\end{align*}
\end{proof}
We recall that  $D(G_0)$ is the 
linear span 
of real and imaginary parts of the functions $u$ of the type $u(t,x) 
=\Phi_j(t) e^{i\langle x,h_j(t)\rangle}$ with 
$\Phi_j\in C_c^1(\R)$ and $h_j\in C^1_b(\R;\R^n)$.

\begin{lemma} 
    \label{DG0dense}
    $D(G_0)$ is dense in $L^{p}(\R^{1+n},\nu)$, for every $p\in [1, 
    +\infty)$.
\end{lemma}
\begin{proof}
Since $\nu$ is a $\sigma$-finite measure on $ \R^{1+n}$ then the 
space of the continuous functions 
with compact support is dense in $L^{p}(\nu)$. Each 
continuous function with compact support $\Phi$ may be approximated 
in the sup 
norm (and hence in  $L^{p}(\nu)$) by a sequence of functions  that 
are 
linear combinations of products $g(t) \varphi(x) $, where both $g$ 
and $\varphi$ 
have compact support. In its turn, each continuous $\varphi$ 
with compact support is the pointwise limit of a sequence of 
exponential functions 
 $\varphi_k $ such that $\|\varphi_{k}\|_{\infty} \leq \|\varphi 
 \|_{\infty}$ for each $k$. The functions $(t,x)\mapsto g(t) 
\varphi_{k}(x) $
  belong to $D(G_{0})$ and approximate $g(t) 
\varphi(x) $ in $L^{p}(\nu)$.  Note 
that 
if $\varphi \in L^{p_{1}}\cap L^{p_{2}}(\R^{1+n},\nu)$, then the 
approximation is simultaneous, i.e. the same sequence approximates 
$\varphi$ both in $L^{p_{1}}(\R^{1+n},\nu)$ and in 
$L^{p_{2}}(\R^{1+n},\nu)$. 
\end{proof}

Since $G_0$ is dissipative in $L^2(\R^{1+n},\nu)$ (see the 
introduction) 
then  it is closable. 
Let us denote the closure of $G_0$ by $G$. Then $G$ is a dissipative, 
densely defined,
closed operator. Moreover, it is the generator of the semigroup
$\cP_\tau$ as the next proposition shows.

\begin{proposition}
	\label{prp:semigroup}
$\mathcal{P}_{\tau}$ is a strongly continuous contraction semigroup 
in $L^2(\R^{1+n}, \nu)$, that leaves  $D(G_{0})$ invariant. Its 
infinitesimal generator is the closure $G$ of $G_{0}$. 
Moreover $\nu$ is an invariant measure for $\mathcal{P}_{\tau}$. 
\end{proposition}
\begin{proof}
The proof is similar to the one in the $T$-periodic case that can be 
found in \cite{DPL06}; we write it here for the reader's convenience. 

Let $u\in D(G_{0})$, $u(t,x) =  \Phi(t)e^{i\langle x,h(t)\rangle}$ 
and fix  
 $\tau >0$. Then  we have
 \begin{equation}
 \label{e4.7bis}
 \begin{array}{l}
 \displaystyle \mathcal{P}_{\tau} 
u(t,x)=\int_{\R^n}\Phi(t+\tau)e^{i\langle 
U(t+\tau,t)x  
 +g(t+\tau,t)+y,h(t+\tau) \rangle}{\mathcal N} _{0,Q(t+\tau,t)}(\ddc 
y)
 \\
 \\
 =\Phi(t+\tau)e^{i\langle g(t+\tau,t),h(t+\tau) \rangle} e^{i\langle  
 U(t+\tau,t)x ,h(t+\tau) \rangle}
 e^{-\frac12\;\langle Q(t+\tau,t) h(t+\tau), h(t+\tau)\rangle}  
 \\
 \\
 :=\Psi_{\tau }(t)e^{i\langle x ,U^{*}(t+\tau,t)h(t+\tau) \rangle}.
 \end{array}
 \end{equation}
Therefore, $\mathcal{P}_{\tau}$ preserves $D(G_{0})$, the semigroup 
law follows easily, 
as well as the strong continuity on $D(G_{0})$.

Let us identify the generator of $\mathcal{P}_{\tau}$  as $G$. 
The domain   $D(G_0)$
is contained in the domain of the infinitesimal generator $L$ of 
$\mathcal{P}_{\tau}$, because  for 
$u = \phi(t)e^{i\langle x,h(t) \rangle}$ we have by \eqref{e4.7bis}
$$\begin{array}{l}
\displaystyle {\frac{d}{d\tau} \mathcal{P}_{\tau}u_{|\tau =0}}
\\
\\
= (\phi'(t) + i\phi(t) \langle x,h'(t) \rangle)e^{i\langle x,h (t) 
\rangle}-
 \bigg(\frac12\;|B^*(t)h(t)|^2+i\langle A(t)x+f(t),h(t) \rangle 
\bigg)u (t,x) 
 \\
 \\
= (G_0u)(t,x).
\end{array} $$
Since $D(G_{0})$ is invariant under $\mathcal{P}_{\tau}$ 
and dense in $L^2(\R^{1+n},\nu)$, then it is a core for  $L$, which 
means that it is dense in 
$D(L)$ for the graph norm. Therefore, $L$ is the closure of $G_0$. 
Since $L=G$ is dissipative, $\cP_\tau$ is a contraction semigroup in
$L^2(\R^{1+n},\nu)$.

The fact that $\nu$ is invariant for $\cP_\tau$ follows easily from \eqref{invariantmess}.
\end{proof}


\section{Smoothing properties of the evolution operator and of the 
evolution semigroups}
\label{secevol}


In this section  we prove estimates for the spatial derivatives of  
 $P_{s,t}\varphi$ with $\varphi\in L^p(\R^{n},\nu_t)$ and for the 
spatial derivatives of 
$\cP_{\tau}u$,    with $u\in L^p(\R^{1+n},\nu)$ or 
$u\in L^p_{\#}(  \R^{1+n},\nu)$. In
order to do so, we first obtain estimates for the
spatial derivatives of the density $\rho_{ U(t,s)x-g(t,s),Q(t,s)}$ of
$\NN_{U(t,s)x-g(t,s),Q(t,s)}$.    
For notational reasons we suppress that $\rho$ depends on $t,s$
and $x$ and we shortly write $\rho(y)$ for $\rho_{
U(t,s)x-g(t,s),Q(t,s)}(y)$.
 
\begin{lemma}
    \label{ndestimate}
    Let $\alpha\in\N_0^n$ with $|\alpha|=k$ and let 
$p,q\in(1,\infty)$ satisfy 
    $1/p+1/q=1$, or $(p,q)=(1,+\infty)$. 
    Then there exists $C>0$ such
    that 
    \begin{equation}        
\|\ND^{-\frac1p} \D^\alpha_x\rho \|_{L^{q}(\R^n,\ddc x)}\leq
C \|Q^{-\frac12}(t,s)\|^k\|U(t,s)\|^k,\quad x\in\R^n,\ t>s.
    \end{equation}
\end{lemma}
\begin{proof}
Since $\D_{x} \rho  = 
\rho \cdot U^{*}(t,s)Q^{-1}(t,s)(y-U(t,s)x+g(t,s))$, 
differentiating further we obtain that  $|\D^\alpha_x \rho|  \leq
\rho \cdot P_{\alpha}(\|\cA_1(t,s,x,y)\|, 
\|\cA_2(t,s,x,y)\|)$, where 
\begin{align*}
  \cA_1(t,s,x,y)&= U(t,s)^{*}Q^{-1}(t,s)(y-U(t,s)x+g(t,s)),& 
s,t\in\R,\
  x,y\in\R^n,\\
  \cA_2(t,s,x,y)&=-U^*(t,s)Q^{-1}(t,s)U(t,s),&s,t\in\R,\ x,y\in\R^n.
\end{align*}
and $P_{\alpha}(\xi, \eta)=\sum_{i+2j=k}\beta_{ij}\xi^i\eta^j$ for 
some
$\beta_{ij}\in\R$.
 
The statement follows now immediately if $p=1$, $q=\infty$. If $p>1$ 
then
\begin{align*}    
 \int\limits_{\R^n}|\ND^{-\frac1p} (y)\D^\alpha_x\rho (y)|^{q}\dd
y    
\leq   C\sum\limits_{i,j\in\N_0,i+2j=k}\,\int\limits_{\R^n}\rho 
(y)    
 \|\cA_1(t,s,x,y) \|^{iq} \|\cA_2(t,s,x,y) \|^{jq}\dd y.
\end{align*}
By the substitution $y=Q(t,s)^{\frac12}\eta+U(t,s)x+g(t,s)$, we obtain
\begin{align*}    
&\int\limits_{\R^n}\rho (y) \|\cA_1(t,s,x,y) \|^{iq} 
\|\cA_2(t,s,x,y) \|^{jq}\dd
y\\
=& (2\pi)^{-\frac n2} \int\limits_{\R^n}\exp(-|\eta|^{2}/2)            \|          
Q^{-\frac12}(t,s)U(t,s)\|^{iq}|\eta|^{iq}\|U^*(t,s)Q^{-1}(t,s)U(t,s)\|^{jq}\dd         
\eta
\\    
\leq&C \|Q^{-\frac12}(t,s)U(t,s) 
\|^{iq}\|U^*(t,s)Q^{-1}(t,s)U(t,s)\|^{jq}\leq    
C\|Q^{-\frac12}(t,s)\|^{iq+2jq}\|U(t,s)\|^{iq+2jq}\\
    =&C\|Q^{-\frac12}(t,s)\|^{kq}\|U(t,s)\|^{kq}.
\end{align*}
Summing up, the proof is complete.
\end{proof}

The next lemma provides estimates for $Q^{-\frac12}(t,s)$.
\begin{lemma}
    \label{qdecay}
    There exist $C,\delta>0$ such that 
    \begin{align*}
        \| Q^{-\frac12}(t,s)\|\leq 
        \left\{
        \begin{array}{rl}
            C(t-s)^{-\frac12},& 0<t-s<\delta,\\
            C,& t-s\geq\delta.
        \end{array}
        \right.
    \end{align*}
\end{lemma}
\begin{proof}
Let $x\in\R^n$. Then, by \eqref{invertb},
\begin{align*}
    \langle Q(t,s)x,x\rangle &=\int\limits_s^t\langle 
U(t,r)B(r)B^*(r)U^*(t,r)x,x\rangle \dd
    r=\int\limits_s^t\|B^*(r)U^*(t,r)x\|^2\dd r\\
    &\geq \mu_0
    \int\limits_s^t\|U^*(t,r)x\|^2\dd r.
\end{align*}
Since $\|U^*(t,r)x-x\|\leq \frac{1}{2}\|x\|$ for $t-r<\delta$ with
some $\delta>0$, independent of $t$, $r$, and $x$, we obtain
\begin{align*}
    \langle Q(t,s)x,x\rangle \geq\frac{\mu_0}{4}(t-s)\|x\|^2,\quad 
0<t-s<\delta.
\end{align*}
Similarly, for $t-s\geq\delta$, we have
$$\langle Q(t,s)x,x\rangle = \int\limits_s^t\|B^*(r)U^*(t,r)x\|^2\dd
    r\geq\int\limits_{t-\delta}^t\|B^*(r)U^*(t,r)x\|^2\dd r \geq 
    \frac{\mu_0}{4}\delta\|x\|^2 . $$
\end{proof}

By \eqref{exp_decay}  we obtain, for $\omega \in (\omega_{0}(U),0)$, 

\begin{equation}
    \label{stimaQ}
    \begin{array}{lll}
    \langle Q(t,s)x,x\rangle& 
=&\displaystyle{\int_{s}^{t}\|B^*(r)U^*(t,r)x\|^2\dd
    r\leq CM\int_{s}^{t}{\rm e}^{ 2\omega(t-r)}\|x\|^2\dd
    r}\\ & \leq&
    \displaystyle{CM\frac{1-{\rm e}^{ 
\omega(t-s)}}{2|\omega|}\|x\|^2,\quad t>s.}
\end{array}
\end{equation}
Hence, $\|Q^{-\frac12}(t,s)\|$ does not decay for $ t-s \to\infty$. 
In 
other
words, the estimate of Lemma~\ref{qdecay} is optimal for $ t-s 
\to\infty$.

Now we are in the position to prove estimates for the spatial
derivatives of   $P_{s,t}\varphi$, for each $\varphi\in 
L^p(\R^{n},\nu_t)$.

\begin{lemma}
    \label{evolutionsysdecay}
    Let $\alpha\in\N_0^n$ and $p\in [1,\infty)$. For
    $\omega>\omega_0(U)$ there exist 
$C,\delta>0$
    such that,
    \[
    \|\D_x^\alpha 
P_{s,t}\|_{\cL(L^p(\R^{n},\nu_t),L^p(\R^{n},\nu_s))}\leq
    \left\{
    \begin{array}{rl}
        C(t-s)^{ |\alpha|/2}{\rm
        e}^{ \omega|\alpha|(t-s)},&0<t-s<\delta,\\
        C{\rm e}^{ \omega|\alpha|(t-s)},&t-s>\delta.
    \end{array}
    \right.
    \]
\end{lemma}
\begin{proof}
Since $C_b(\R^n)$ is dense in $L^p(\R^{n},\nu_t)$ for any $t\in\R$, 
it is enough to estimate $\|\D_x^\alpha 
P_{s,t}\|_{L^p(\R^{n},\nu_s)}$ 
for $\varphi\in C_b(\R^n)$. H\"older's inequality and 
Lemma~\ref{ndestimate} yield
\begin{align*}
    |\D_x^\alpha    
P_{s,t}\varphi(x)|^p=& 
\;\bigg|\D_x^\alpha\int\limits_{\R^n}\varphi(y+g(t,s))\rho (y)\dd
y\bigg|^p\\
 \leq & \;
 \int\limits_{\R^n}|\varphi(y+g(t,s))|^p\rho (y)\dd y \,   
\|\ND^{-\frac1p} \D^\alpha_x\rho \|_{L^{q}(\R^n,\ddc x)}^p\\    
=& \; 
 \left(P_{s,t}|\varphi|^p\right)(x)    
\|\ND^{-\frac1p} \D^\alpha_x\rho \|_{L^{q}(\R^n,\ddc x)}^p\\
\leq& \;C
 \left(P_{s,t}|\varphi|^p\right)(x)
 \|Q^{-\frac12}(t,s)\|^k\|U(t,s)\|^k,
\end{align*}
where $1/p+1/q=1$.
Hence,  it follows from \eqref{invariantmess} that 
\begin{align*}
 \|\D_x^\alpha 
P_{s,t}\varphi\|_{L^p(\nu_t)}^p\leq & \;C \|Q^{-\frac12}(t,s)\|^k\|U(t,s)\|^k
    \int\limits_{\R^n}P_{s,t}|\varphi|^p(x)\nu_t(\ddc x)\\
= & \;C \|Q^{-\frac12}(t,s)\|^k\|U(t,s)\|^k\int\limits_{\R^n}
|\varphi|^p(x)\nu_s(\ddc x)\\
= & 
\;C \|Q^{-\frac12}(t,s)\|^k\|U(t,s)\|^k\|\varphi\|_{L^p(\R^n,\nu_s)}^p.
\end{align*}
Here, we have used that $|\varphi|^p\in C_b(\R^n)$.
Now, Lemma~\ref{qdecay} and \eqref{exp_decay} yield the assertion.
\end{proof}

Thanks to the representation \eqref{rep}, the smoothing properties of 
the evolution operator $P_{s,t}$, given in
Lemma~\ref{evolutionsysdecay}, 
yield smoothing properties of the semigroups $(\cP_\tau)_{\tau\geq0}$ 
and $(\cP_\tau^\#)_{\tau\geq0}$.

\begin{lemma}
    \label{decaysemigroup}
Let $\alpha\in\N_0^n$ and $p\in [1,\infty)$. For $\omega>\omega_0(U)$ 
there exist  
$C,\delta>0$, such that
    \begin{align*}
	\|\D_x^\alpha\cP_\tau\|_{\cL(L^p(\R^{1+n},\nu))}\leq
	\left\{
	\begin{array}{rl}
	    C\tau^{-\frac{|\alpha|}{2}}{\rm
	    e}^{\omega|\alpha|\tau}&,\quad
	     0<t-s<\delta,\\
	     C{\rm e}^{\omega|\alpha|\tau}&,\quad
	     t-s\geq\delta.
	\end{array}
	\right.
    \end{align*}
\end{lemma}
\begin{proof}
By Lemma~\ref{evolutionsysdecay}, there exist $C,\delta>0$, such that 
\begin{align*}
    \|\D^\alpha_x\cP_\tau    
u\|^p_{L^p(\R^{1+n},\nu)}=&\int\limits_\R\int\limits_{\R^n}|\D_x^\alpha\cP_\tau    
u(t,x)|^p\nu_t(\ddc x)\dd 
t=\int\limits_\R\int\limits_{\R^n}|\D_x^\alpha
    P_{t,t+\tau}u(t+\tau ,x)|^p \nu_t(\ddc x)\dd
    t\\\leq&CK(\tau)^p\int\limits_\R \int\limits_{\R^n}|
    u(t+\tau,x)|^p \nu_{t+\tau}(\ddc x)\dd t,\quad u\in
    L^p(\R^{1+n},\nu),
\end{align*}
where 
\begin{align*}
 K(\tau):=
 \left\{
 \begin{array}{rl}
   \tau^{-\frac{|\alpha|}{2}}{\rm e}^{\omega|\alpha|\tau}&,\quad 
0<t-s<\delta,\\
  {\rm e}^{\omega|\alpha|\tau}&,\quad t-s\geq\delta.
 \end{array}
 \right.
\end{align*}
Now, the lemma follows from the substitution $s=t+\tau$. Indeed, for 
each $v \in L^p(\R^{1+n},\nu)$ we have
\begin{align*}
    \int\limits_\R \int\limits_{\R^n}|v (t+\tau,x)|^p
    \nu_{t+\tau}(\ddc x)\dd t =& \int\limits_\R
    \int\limits_{\R^n}|v (s,x)|^p\nu_{s}(\ddc x)\dd s=
    \|v \|_{L^p(\R^{1+n},\nu)}^p  .  
\end{align*}
\end{proof}

The proof of Lemma~\ref{decaysemigroup} can be easily carried
over to the $T$-periodic case and it is therefore omitted.

\begin{lemma}
    \label{sharpdecaysemigroup}
    Let $\alpha\in\N_0^n$ and $p\in [1,\infty)$. For each
    $\omega>\omega_0(U)$ there exist $C,\delta>0$, such	that
\begin{align*}
\|\D_x^\alpha\cP _\tau\|_{\cL(L^p_\#(\R^{1+n},\nu))}\leq
	\left\{
	\begin{array}{rl}
	    C\tau^{-\frac{|\alpha|}{2}}{\rm
	    e}^{\omega|\alpha|\tau}&,\quad
	     0<t-s<\delta,\\
	     C{\rm e}^{\omega|\alpha|\tau}&,\quad
	     t-s\geq\delta.
	\end{array}
	\right.
    \end{align*}
\end{lemma}


\section{The spaces $H^{k,s}(\R^{1+n},\nu)$ and 
$H_\#^{k,s}(\R^{1+n},\nu)$}
\label{secinterpol}


In this section we show some properties of the spaces
$H^{k,s}(\R^{1+n},\nu)$ and
$H_\#^{k,s}( \R^{1+n},\nu)$ defined in the introduction. 

We recall that for any Gaussian measure $\NN_{m,Q}$ and for any 
$k\in\N_0$, the space $H^k(\R^n,$ $ \NN_{m,Q})$ is defined as 
\begin{align*}
    H^k(\R^n,\NN_{m,Q}):=
    \left\{ f\in L^2(\R^n,\NN_{m,Q}):\exists\D_\beta f\in
    L^2(\R^n,\NN_{m,Q}),\ |\beta|\leq k \right\}.
\end{align*}
Then, $H^k(\R^n,\NN_{m,Q})$ equipped with its natural norm is a 
Hilbert space.

We first show that $D(G_0)$ is dense in $H^{1,2}(\R^{1+n},\nu)$.

\begin{lemma} \label{dense1}
    $D(G_0)$ is dense in  
    $H^{1,2}(\R^{1+n},\nu)$.
\end{lemma}
\begin{proof}
Let $v \in C^\infty_c(\R^{1+n})$. We choose $R_0>0$ and $S>0$ 
such that
$\supp v \subset \Omega_{S,R_0}$, where 
$\Omega_{S,R_0}:=(-S,S)\times(-R_0,R_0)^n$. 
For $R>R_0$ and $k\in \Z^{n}$ we set
\begin{align*}
    v _{R,l}(t,x)=\sum\limits_{|k|=0}^la_{R,k}(t){\rm
    e}^{i\frac R\pi\langle k,x\rangle },
\end{align*}
where $a_{R,k}(t):=\frac1{(2R)^n}\int_{(-R,R)^{n}} v (t,x){\rm 
e}^{-i\frac
R\pi\langle k,x\rangle }\dd x$.

We will show that for each
$\varepsilon >0$ there 
exists
$R>R_0$ and $l\in\N$ such that
\begin{align*}
\|v -v _{R,l}\|_{H^{1,2}(\R^{1+n},\nu)}\leq \varepsilon.
\end{align*}

Since $v$ is smooth and compactly supported, there exists $K>0$ such 
that for each $R\geq R_{0}$ and $k\in \Z^{n}$ we have
\begin{align*}    
 \| \partial_t v_{R,l} \|^{2}_{L^\infty(\R^{1+n})}+ 
\sum\limits_{|\alpha|=0}^2\|D^\alpha_x v_{R,l} 
\|^{2}_{L^\infty(\R^{1+n})} \leq K.
\end{align*}
Let us fix $R\geq R_{0}$ such that $K \nu( (-S,S)\times\R^n\setminus \Omega_{S,R} ) \leq \varepsilon 
/2$. As $l\to \infty$, the sequences $(v_{R,l})$, $(\partial_{t} 
v_{R,l})$, $(\D_x^\alpha v _{R,l})$ converge uniformly on $ 
\Omega_{S,R}$ 
to $v$, $ \partial_{t}
v$, $\D_x^\alpha v$, respectively. Therefore, there exists $l\in\N$ such that
\begin{align*}    
\sum\limits_{|\alpha|=0}^2\|\D_x^\alpha v _{R,l}-\D_x^\alpha v 
\|^2_{L^\infty(
\Omega_{S,R})}+\|\partial_tv _{R,l}-\partial_t v \|^2_{L^\infty(
\Omega_{S,R})}\leq \frac\varepsilon{2\nu(\Omega_{S,R}) } .
\end{align*}
For such $l$ we have 
\begin{align*}    
\|v -v _{R,l}\|^2_{H^{1,2}(\R^{1+n},\nu)}=&\|v -v _{R,l}\|^2_{H^{1,2}(
\Omega_{S,R},\nu)}+\|v _{R,l}\|^2_{H^{1,2}(
    \Omega_{S,R}^c,\nu)}
		\leq \frac\varepsilon2
  +\frac\varepsilon2=\varepsilon.
    \end{align*}
    Since $C_c^\infty(\R^{1+n})$ is dense in
    $H^{1,2}(\R^{1+n},\nu)$, the proof is complete.
\end{proof}

The corresponding result for $T$-periodic spaces reads as
follows.

\begin{lemma} \label{sharpdense1}
    $D(G_0^\#)$ is dense in $L_\#^{ 2}(\R^{1+n},\nu)$ and in 
    $H_\#^{1,2}(\R^{1+n},\nu)$.
\end{lemma}
\begin{proof}
Let 
\begin{align*}
v \in C^\infty_{c,\#}:=\big\{u\in
C^\infty(\R^{1+n})&:u(t,x)=u(T+t,x)\mbox{ for all }
t\in\R, x\in\R^n\mbox{, and }\\&\supp
u(t,\cdot)\subset(-R_0,R_0)^n\mbox{ for all }
t\in\R\mbox{ and some }R_0>0\big\}.
\end{align*}
We define
\begin{align*}   
K:=\|\partial_tv \|_{L^\infty(\R^{1+n})} + 
\sum\limits_{|\alpha|=0}^2\|D^\alpha_xv \|_{L^\infty(\R^{1+n})}.
\end{align*}
and, for $R>R_0$, we set
\begin{align*}
    v _{R,l}(t,x)=\sum\limits_{|k|=-l}^la_{R,k}(t){\rm
    e}^{i\frac R\pi\langle k,x\rangle },
\end{align*}
where $a_{R,k}(t):=\frac1{(2R)^n}\int_{-R}^Rv (t,x){\rm 
e}^{i\frac
R\pi\langle k,x\rangle }\dd x$. Clearly, $a_{R,k}$ is $T$-periodic. 
As 
in the 
proof of
Lemma~\ref{dense1}, it follows that for $\varepsilon\in(0,K)$ there 
exists
$R>R_0$ and $l\in\N$ such that
\begin{align*}
\|v -v _{R,l}\|_{H_\#^{1,2}(\R^{1+n},\nu)}\leq 
\varepsilon.
\end{align*}
Since $C_{c,\#}^\infty((0,T)\times\R^{n})$ is dense in
$H_\#^{1,2}((0,T)\times\R^{n},\nu)$, the proof is complete.
\end{proof} 

Let $t_0\in\R$.
In the following we denote the product measure of the one
dimensional Lebesgue measure and $\nu_{t_0}$ by $\ddc 
t\times\nu_{t_0}$.

\begin{lemma}
\label{lem:isom}
\begin{enumerate}
\item
There exists an isomorphism
\begin{align*}
\cT: L^2(\R^{1+n},\nu)\to L^2(\R^{1+n},\ddc t \times \NN_{0,1}),
\end{align*}
such that, for $k=0$, $s\in\N_0$ and
for $k=1$, $s=2$,
$\cT|_{H^{k,s}(\R^n,\nu)}$ is an isomorphism from
$H^{k,s}(\R^{1+n},\nu)$ onto $H^{k,s}(\R^{1+n},\ddc t \times
\NN_{0,1})$.
\item
Let $t_0\in\R$. Then, there exists an isomorphism
\begin{align*}
\cT_{t_0}: L^2(\R^{1+n},\nu)\to L^2(\R^{1+n},\ddc
t\times\nu_{t_0}),
\end{align*}
such that, for $k=0$, $s\in\N_0$ and for
$k=1$, $s=2$,
$\cT_{t_0}|_{H^{k,s}(\R^n,\nu)}$ is an isomorphism from
$H^{k,s}(\R^{1+n},\nu)$ onto $H^{k,s}(\R^{1+n},\ddc
        t\times\nu_{t_0})$.
\end{enumerate}
\end{lemma}
\begin{proof}
For $t\in\R$ and $x\in\R^n$ define 
\begin{align*}
\left( \cT u \right)(t,x) :=
u\left(t,Q^{\frac12}(t,-\infty)
x+g(t,-\infty)\right).
\end{align*}
By substitution, we obtain
\begin{align*}
\|u\|_{L^2(\R^{1+n},\nu)}=\|\cT u\|_{L^2(\R^{1+n},\ddc t \times 
\NN_{0,1})},\quad 
u\in L^2(\R^{1+n},\nu).
\end{align*}
Moreover, for $l\in\N$ and for arbitrary integers $\alpha_{1}, \ldots 
\alpha_{l} \in \{1, \ldots ,n\}$  we have
\begin{align*}
   \frac{\partial  \cT u}{\partial x_{\alpha_{1}}\ldots \partial 
   x_{\alpha_{l}}}\, (x) = 
	 \sum_{\beta_i\in\{1,\dots,n\}}\frac {\partial 
u}{\partial_{\beta_{1}}\ldots \partial_{\beta_{l}}}(Q^{ 
\frac12}(t,-\infty)x
+g(t,-\infty))\prod_{i=1}^{l}(Q^{ 
\frac12}(t,-\infty)_{\beta_{i}\alpha_{i}} 
\end{align*}
so that, for $u\in H^{0,k}(\R^{1+n},\nu)$,  
\begin{align*}
\|\,|\D_x^k \cT u|\,\|_{L^2(\R^{1+n},\ddc t \times \NN_{0,1})} \leq 
&\|Q^{ \frac12}(t,-\infty)\|^{k}
\| \,|\D^k_xu|\, \|_{L^2(\R^{1+n},\nu)}.
\end{align*}
and similarly, for $ u\in  H^{0,k}(\R^{1+n},\ddc t \times \NN_{0,1})$,
\begin{align*}
\|\,|\D^k_x\cT^{-1}u|\,\|_{L^2(\R^{1+n},\nu)} \leq &
 \|Q^{-\frac 12}(t,-\infty)\|^{k}
 \|\,|\D_x^k u|\,\|_{L^2(\R^{1+n},\ddc t \times \NN_{0,1})}.
\end{align*}
Since the norms $\|Q^{\frac 12}(t,-\infty)\| $ and $\|Q^{-\frac 
12}(t,-\infty)\|$ are bounded by a constant independent of $t$, then 
$\cT$ is an isomorphism from 
$H^{0,k}(\R^{1+n},\nu)$ to $H^{0,k}(\R^{1+n},\ddc t \times 
\NN_{0,1})$ for $k\in\N_0$.

For $\varphi\in H^{1,2}( \R^{1+n},\nu)$  we have
\begin{align*}
\partial_t(\cT\varphi)(t,x)
=&\partial_t\varphi(t,Q^\frac12(t,-\infty)x-g(t,-\infty))
=(\partial_t\varphi) (t,Q^\frac12(t,-\infty)x-g(t,-\infty))\\
&+\langle (\nabla_{x}\varphi) (t,Q^\frac12(t,-\infty)x-g(t,-\infty)), 
\partial_tQ^\frac12(t,-\infty)x-\partial_t
g(t,-\infty)\rangle.
\end{align*}
Clearly, $\partial_t 
g(t,-\infty)=f(t)+\int_{-\infty}^tA(t)U(t,r)f(r)\dd r$ is
uniformly bounded for $t\in\R$. Moreover, the representation
\begin{align*}
Q^\frac12(t,-\infty)=\frac1{2\pi i}\int\limits_\Gamma
\lambda^\frac12 R(\lambda,Q(t,-\infty))\dd \lambda
\end{align*}
for a suitable path $\Gamma$ yields
\begin{align*}
\|\partial_t Q^\frac12(t,-\infty)\|\leq C,\quad t\in\R,
\end{align*}
thanks to the uniform boundedness of $\partial_t
R(\lambda,Q(t,-\infty))$ for $\lambda\in\Gamma$ and $t\in\R$. The 
latter follows from the boundedness of $\|A(t)\|$, $\|B(t)\|$, and 
$\|Q(t, -\infty)\|$, see \eqref{stimaQ}. 
Moreover, an easy computation (see e.g. \cite[Lemma 2.1]{Lun97T}, or  \cite[Lemma 2.3]{MPRS02})
shows that  there exists $C>0$ such 
that 
for   any $M \in\cL(\R^n)$ and $t\in\R$, $\psi \in H^2(\R^n,\nu_t)$
\begin{equation}
    \|\langle M \cdot,\D_x \psi\rangle \|_{L^2(\R^n,\nu_t)}\leq    
C \|M \|  \cdot \|\psi\|_{H^2(\R^n,\nu_t)}.
    \label{unibound}
\end{equation}
Therefore, it follows that 
\begin{align*}
\|\cT\varphi\|_{H^{1,2}( \R^{1+n},\ddc t \times \NN_{0,1})}
\leq C
\|\varphi\|_{H^{1,2}(\R^{1+n},\nu)}, \quad \varphi\in
H^{1,2}(\R^{1+n},\nu).
\end{align*}
Similarly, we obtain
\begin{align*}
\|\cT^{-1}\varphi\|_{H^{1,2}(\R^{1+n},\nu)}
\leq C
\|\varphi\|_{H^{1,2}( \R^{1+n},\ddc t \times \NN_{0,1})},
\quad \varphi \in
H^{1,2}(\R^{1+n},\ddc t \times \NN_{0,1}).
\end{align*}
This proves (a). 

Setting
\begin{align*}
\left( \cT_{t_0} u \right)(t,x) :=
u\left(t,Q^{\frac12}(t,-\infty)Q^{-\frac12}(t_0,-\infty)
\left(x-g(t_0,-\infty)\right)+g(t,-\infty)\right),
\end{align*}
assertion (b) follows as above.
\end{proof}
There is a corresponding result for T-periodic spaces as well. Since 
the proof of the
following lemma is similar to the proof of Lemma~\ref{lem:isom} it is
omitted.
\begin{lemma}
\label{lem:isom2}
\begin{enumerate}
\item
There exists an isomorphism
\begin{align*}
\cT_\#: L_\#^2( \R^{1+n},\nu)\to L^2_\#( \R^{1+n},\ddc t \times 
\NN_{0,1}),
\end{align*}
such that, for $k=0$ and $s\in\N_0$ and
$k=1$ and $s=2$,
$\cT_\#|_{H_\#^{k,s}( \R^{1+n},\nu)}$ is an isomorphism from
$H^{k,s}_\#( \R^{1+n},\nu)$ onto $H_\#^{k,s}( \R^{1+n},\ddc t \times
\NN_{0,1})$. 
\item
Let $t_0\in\R$. Then, there exists an isomorphism
\begin{align*}
\cT_{\#,t_0}: L^2_\#( \R^{1+n},\nu)\to L_\#^2(
\R^{1+n},\ddc
t\times\nu_{t_0}),
\end{align*}
such that, for $k=0$ and $s\in\N_0$ and
$k=1$ and $s=2$, 
$\cT_{\#,t_0}|_{H^{k,s}_\#( \R^{1+n},\nu)}$ is an isomorphism from
$H_\#^{k,s}( \R^{1+n},\nu)$ onto $H_\#^{k,s}(
\R^{1+n},\ddc
        t\times\nu_{t_0})$.
\end{enumerate}
\end{lemma}
Next, we give a characterization of some real interpolation spaces 
between $L^2(\R^{1+n},\nu)$ and $H^{0,s}(\R^{1+n},\nu)$.

\begin{proposition}
    \label{interpolspaces}
    Let $r,s\in\N$ with $0<r<s$. Then, 
    \begin{align*}
        \left(
        L^2(\R^{1+n},\nu),H^{0,s}(\R^{1+n},\nu)
        \right)_{\frac rs,2}=H^{0,r}(\R^{1+n},\nu).
    \end{align*}
\end{proposition}
\begin{proof}
Let $\cT$ be as in Lemma~\ref{lem:isom}.
Then,  
\begin{align*}
&    \left( L^2(\R^{1+n},\nu), H^{0,s}(\R^{1+n},\nu)\right)_{\frac 
rs,2}\\
& = \left(\cT^{-1}L^2(\R^{1+n},\ddc t \times 
\NN_{0,1}),\cT^{-1}H^{0,s}(\R^{1+n},\ddc t \times 
\NN_{0,1})\right)_{\frac rs,2}
    \\
& = \cT^{-1}\left(L^2(\R^{1+n},\ddc t \times 
\NN_{0,1}),H^{0,s}(\R^{1+n},\ddc t \times \NN_{0,1})\right)_{\frac 
rs,2}.
\end{align*}
On the other hand, by  \cite[Theorem 1.18.4]{Tri78},
\begin{align*}
&	\left(L^2(\R^{1+n},\ddc t \times \NN_{0,1}),H^{0,s}(\R^{1+n},\ddc t
	\times \NN_{0,1})\right)_{\frac rs,2}\\
&= 
\left(L^2(\R,L^2(\R^{n},\NN_{0,1})),L^2(\R,H^s(\R^n,\NN_{0,1}))\right)_{\frac 
rs,2}
\\    
&= L^2\left(\R,(L^2(\R^n,\NN_{0,1}),H^s(\R^n,\NN_{0,1}))_{\frac
rs,2}\right).
\end{align*}    
The real interpolation spaces 
between $L^{2}$ and $H^{s}$ spaces with respect to the Gaussian 
measure $\NN_{0,1}$ are  known, see e.g. 
\cite[Proposition 4]{FL06}. More precisely, we have 
$$H^l(\R^n,\NN_{0,1})=\left(L^2(\R^n,\NN_{0,1}),H^k(\R^n,\NN_{0,1})\right)_{l/k,2},$$
for each $k,l\in\N$ with $0<l<k$, and, therefore, we get 
\begin{align*}
\left( L^2(\R^{1+n},\nu), H^{0,s}(\R^{1+n},\nu) \right)_{\frac rs,2}
    =&\cT^{-1}L^2\left(\R,H^{r}(\R^n,\NN_{0,1})\right)
=\cT^{-1}H^{0,r}(\R^{1+n},\ddc t \times \NN_{0,1})\\
    =&H^{0,r}(\R^{1+n},\nu).
\end{align*}
\end{proof}

Again, there is a corresponding result for $T$-periodic spaces. Since
the proof of Proposition~\ref{interpolspaces} carries over to the 
$T$-periodic
case with only minor modifications, the proof of the next lemma is 
omitted.

\begin{proposition}
    \label{sharpinterpolspaces}
    Let $r,s\in\N$ with $0<r<s$. Then, 
    \begin{align*}
        \left(        
L^2_\#(\R^{1+n},\nu),H_\#^{0,s}(\R^{1+n},\nu)        
\right)_{\frac rs,2}=H_\#^{0,r}(\R^{1+n},\nu).
    \end{align*}
\end{proposition}


\section{The domains of $G$ and of $G_{\#}$}


The proofs of Theorem~\ref{main1} and Theorem~\ref{sharpmain1} 
are based on the following abstract interpolation result (see
\cite[Theorem 2.5]{Lun99}).

\begin{proposition}
    \label{domainl}
    Let $T(t)$ be a semigroup on some Banach space $X$ with generator 
$L:D(L)\to X$.
    Assume that there exists a Banach space $E\subset X$ and
    $m\in\N$, $0<\beta<1$, $\omega\in\R$, $C>0$ such that
    \begin{align*}
        \|T(t)\|_{\cL(X,E)}\leq C{\rm e}^{\omega
        t}t^{-m\beta},\quad t>0,
    \end{align*}
    and for every $x\in X$, $t\mapsto T(t)x$ is measurable with
    values in $E$. Then $E\in J_\beta(X,D(L^m))$, so that
    $(X,D(L^m))_{\theta,p}\subset(X,E)_{\theta/\beta,p}$, for every
    $\theta\in(0,\beta)$, $p\in [1,\infty]$.
\end{proposition}
Indeed, we apply this proposition taking  $X = L^2(\R^{1+n},\nu)$, 
$T(t) = \cP_t$ and $ E=H^{0,k}(\R^{1+n},\nu)$, or  
$X=L^{2}_\#( \R^{1+n},\nu)$, $T(t) = \cP _t$ and $E=H^{0,k}_\#( 
\R^{1+n},\nu)$ in the periodic case. 
With these choices we have  optimal blow-up estimates for the norms  
$\|T(t)\|_{\cL(X,E)}$  as $t\to 0$, given by Lemmas 
\ref{decaysemigroup} 
and \ref{sharpdecaysemigroup}. Since the real interpolation
spaces between $X$ and $E$ have been characterized 
in Propositions~\ref{interpolspaces} and \ref{sharpinterpolspaces}, 
the main
step of the proof of Theorem~\ref{main1} and Theorem~\ref{sharpmain1}
follows.

\subsection{Proof of Theorem~\ref{main1}}
\label{secmain1}
We first prove the continuous embedding
\begin{align}
    D(G)\subset H^{1,2}(\R^{1+n},\nu).
    \label{step1}
\end{align}
We use Proposition \ref{domainl} with $X=L^2(\R^{1+n},\nu)$ and 
$E=H^{0,4}(\R^{1+n},\nu)$.  
From  Lemma~\ref{decaysemigroup} it follows that there are $C$, 
$\omega$, such that 
$$ \| {\mathcal P}_{\tau}\|_{{\mathcal L}(X,E)} \leq Ce^{\omega 
\tau }\tau^{-2}, \quad \tau >0.$$
Choosing  $m=4$, $\theta=\frac14$ and $\beta=\frac12$, Proposition 
\ref{domainl} yields
\begin{align*}
    (L^2(\R^{1+n},\nu),D(G-I)^4 )_{\frac14,2}\subset \left(
    L^2(\R^{1+n},\nu),H^{0,4}(\R^{1+n},\nu)
    \right)_{\frac12,2}. 
\end{align*}
Since
$L^2(\R^{1+n},\nu)$ is a Hilbert space and $G-I$ is dissipative and 
invertible,
\cite[Theorem 5]{Kat62} and \cite[Theorem 1.15.3]{Tri78} yield
$(L^2(\R^{1+n},\nu), D(G-I)^4)_{\frac14,2}=D(G-I)=D(G)$.
Therefore, Proposition~\ref{interpolspaces} implies
\begin{align}
D(G)\subset H^{0,2}(\R^{1+n},\nu).
 \label{step10}
\end{align}

Now, let $u\in D(G_0)$ and set 
\begin{align*}
\psi (t,x):= \left(G u\right)(t,x)={\frac{\partial}{\partial
t}}u(t,x)+\cK(t)u(t,x),\quad t\in\R,\ x\in\R^n.
\end{align*}
It follows from \eqref{unibound}, that there exists $C_{1}>0$, 
independent of $u$, such that
\begin{align*}
    \|\cK(t)u(t)\|_{L^2(\R^n,\nu_t)}\leq
    C_{1}\|u(t, \cdot)\|_{H^2(\R^n,\nu_t)},\quad
    t\in\R.
\end{align*}
Moreover by \eqref{step10}, 
$\|\cK(\cdot)u\|_{L^2(\R^{1+n},\nu)}\leq 
C_{2}(\|u\|_{L^2(\R^{1+n},\nu)}+\|G u\|_{L^2(\R^{1+n},\nu)})$,
where $C_{2}>0$ is independent of $u$.
Writing $\frac{\partial}{\partial t}u=\psi -\cL(\cdot)u$, we obtain
\begin{align*}
    \|u_{t}\|_{L^{2}(\R^n,\nu)}\leq
  \|\psi \|_{L^2(\R^{1+n},\nu)}+\|\cL(\cdot)u\|_{L^2(\R^{1+n},\nu)}
\\
\leq \|Gu\|_{L^2(\R^{1+n},\nu)} + C_{1}C_{2}\left( \|u \|_{L^{2}(\R^n,\nu)}
+ \|Gu\|_{L^2(\R^{1+n},\nu)}
  \right)=2C \|G u\|_{L^2(\R^{1+n},\nu)}.
\end{align*}
Putting together this estimate and \eqref{step10} we get 
$$\|u\|_{H^{1,2}(\R^{1+n},\nu)} \leq C_{3}\left( \|u \|_{L^{2}(\R^n,\nu)}
+ \|Gu\|_{L^2(\R^{1+n},\nu)}\right)$$
with $C_{3}$ independent of $u$. 
Since $D(G_0)$ is a core of $D(G)$,  
the proof of \eqref{step1} is complete. Moreover, since, by
Lemma~\ref{dense1}, $D(G_0)$ is
dense in $H^{1,2}(\R^{1+n},\nu)$, we have
$D(G)= H^{1,2}(\R^{1+n},\nu)$.

Now let us prove the second equality. The inclusion $``\subset"$ is 
obvious. Let 
\begin{align*}u\in \left\{u\in
H^{1,2}_\mathrm{loc}(\R^{1+n})\cap L^2(\R^{1+n},\nu):\cG
u\in L^2(\R^{1+n},\nu)\right\}, 
\end{align*}
fix $\lambda >0$ and set $\psi :=\lambda u-\cG u$. Then 
$v:=u-R(\lambda,G)\psi$
satisfies $\lambda v-\cG v=0$. We will prove that
$v\equiv0$, and hence $u\in D(G)$, provided $\lambda$ is large enough.

In order to do so, let $\varphi\in C_c^\infty(\R^{1+n})$ be such that 
$\varphi(\cdot, \cdot) \equiv1$ on $[-1,1]\times B(0,1)$ and
$\varphi(\cdot, \cdot)\equiv 0$ 
outside $[-2,2]\times B(0,2)$. 
Then
$\varphi_k(t,x):=\varphi(t/k,x/k)$ satisfies $\varphi_k(t,x)\to 1$ for
$k\to\infty$ and $|\D_x\varphi_k(t,x)|\leq 
\|\,|\D_x\varphi|\,\|_{\infty}$ for $t\in\R$ and
$x\in\R^n$.

A direct calculation yields
\begin{align}
\cG(gh)=g\cG h+h\cG g+\langle B^*\D_x g,B^*D_{x}h\rangle ,\quad
g,h\in H^{1,2}_\mathrm{loc}(\R^{1+n}).\label{formula1}
\end{align}
Hence, we obtain
\begin{align*}
0 &=\int\limits_{\R^{1+n}}(\lambda v-\cG v)\varphi_k^2v\dd\nu
=\lambda\|\varphi_kv\|^2_{L^{2}(\R^{1+n},\nu)}-\int\limits_{\R^{1+n}}\varphi_k\cG 
v\,\varphi_k
v\dd\nu
\\    
&=\lambda\|\varphi_kv\|^2_{L^{2}(\R^{1+n},\nu)}-
\int\limits_{\R^{1+n}}\cG(\varphi_kv)\varphi_kv
\dd\nu+\int\limits_{\R^{1+n}}(\cG\varphi_k)v\varphi_kv\dd\nu\\
&\quad+\int\limits_{\R^{1+n}}\langle B^*\D_x\varphi_k,B^*\D_xv\rangle 
\varphi_kv\dd\nu.
\end{align*}
Since $\int_{\R^{1+n}}\cG g\dd \nu=0$ for each $g\in D(G_0)$, it 
follows from
\eqref{formula1} that
\begin{align}\label{formula2}
    \int\limits_{\R^{1+n}}g\,\cG h\dd\nu+\int\limits_{\R^{1+n}}h\,\cG
    g\dd\nu+\int\limits_{\R^{1+n}}\langle B^*\D_x g,B^*\D_x
    h\rangle \dd\nu=\int\limits_{\R^{1+n}}G(gh)\dd\nu=0
\end{align}
for $g,h\in D(G_0)$. Note that \eqref{formula2} also holds for $g,h\in
D(G)$, since $D(G_0)$ is a core of $D(G)$.
In particular, since $\varphi_k v\in D(G)$, we obtain 
\begin{align*}
 \int\limits_{\R^{1+n}}\cG(\varphi_kv)\varphi_kv\dd\nu
 =&-\frac12\int\limits_{\R^{1+n}}\langle B^*\D_x(\varphi_kv),B^*\D_x
(\varphi_kv)\rangle \dd\nu 
\\ 
 =& -   \frac12\int\limits_{\R^{1+n}}\langle 
B^*\D_xv,B^*\D_x          
v\rangle \varphi_k^2\dd\nu-\int\limits_{\R^{1+n}}\langle 
B^*\D_x\varphi_k,B^*\D_x          
v\rangle v\varphi_k\dd\nu 
\\ 
 &-\frac12\int\limits_{\R^{1+n}}\langle B^*\D_x\varphi_k,B^*\D_x
\varphi_k \rangle v^2\dd\nu 
\\ 
 \leq&  -\frac12\|\varphi_k|B^*\D_xv|\,\|^2_{L^{2}(\R^{1+n},\nu)}
+\frac14\|\varphi_k|B^*\D_xv|\,\|^2_{L^{2}(\R^{1+n},\nu)}\\
&+\frac32\|v|B^*\D_x\varphi_k|\,\|^2_{L^{2}(\R^{1+n},\nu)}.
\end{align*}
Hence,
\begin{align*}    
0&\geq\lambda\|\varphi_kv\|^2_{L^{2}(\R^{1+n},\nu)}
+\frac14\|\varphi_k|B^*\D_xv|\,\|^2_{L^{2}(\R^{1+n},\nu)}
-\frac32\|v|B^*\D_x\varphi_k|\,\|^2_{L^{2}(\R^{1+n},\nu)}\\
&\quad-C\|v\|^2_{L^{2}(\R^{1+n},\nu)}-\frac18\|\varphi_k|B^*\D_xv|\,\|^2_{L^{2}(\R^{1+n},\nu)}
-2\|v|B^*\D_x\varphi_k|\,\|^2_{L^{2}(\R^{1+n},\nu)}\\
&\geq    
\lambda\|\varphi_kv\|^2_{L^{2}(\R^{1+n},\nu)}-C_{1}\|v\|^2_{L^{2}(\R^{1+n},\nu)}.
\end{align*}
Letting $k\to\infty$, we obtain
$0\geq(\lambda-C_{1})\|v\|^2_{L^{2}(\R^{1+n},\nu)}$, which implies
$v\equiv0$ provided $\lambda$ is large enough.
\qed \\[2mm]

\subsection{Proof of Theorem~\ref{sharpmain1}}
\label{secmain2}
The proof of Theorem  \ref{sharpmain1} is the same of Theorem 
\ref{main1}, with the space $E= H^{0,4}_\#( \R^{1+n},\nu)$ 
instead of $H^{0,4}(\R^{1+n},\nu)$
and using Lemma \ref{sharpdecaysemigroup} instead of Lemma 
\ref{decaysemigroup} in the first part, and functions $\varphi_{k}$ 
depending only on $x$ in the second part: $\varphi_{k}(x) = 
\varphi(x/k)$, 
with $\varphi \in C^{\infty}_{c}(\R^{n})$ such that $\varphi\equiv 1$ 
on $B(0,1)$ and $\varphi\equiv 0$ outside $B(0,2)$. We omit it. 

\vspace{2mm}

The characterization of $D(G_\#)$ given by Theorem~\ref{sharpmain1} 
implies also that $D(G_\#)$ is compactly embedded in $L^2_\#(
\R^{1+n},\nu)$, through the next Proposition.

\begin{proposition}
\label{compact_emb}
$H^{1,2}_\#( \R^{1+n},\nu)$ is compactly embedded in 
$L^2_\#(\R^{1+n},\nu)$.
\end{proposition}
\begin{proof}
Let $\cT_{\#}$ be as in Lemma~\ref{lem:isom2}.
Writing $H^{1,2}_\#( \R^{1+n},\nu)=\cT_{\#}^{-1}H^{1,2}_\#(
\R^{1+n},\ddc t \times \NN_{0,1})$, it suffices to show that 
$H^{1,2}_\#(
\R^{1+n},\ddc t \times \NN_{0,1})$ is compactly embedded in
$L^2_\#(\R^{1+n},\ddc t \times \NN_{0,1} )$.

Let $u\in B$, where $B$ denotes the unit ball in
$H^{1,2}_\#(\R^{1+n},\ddc t \times \NN_{0,1})$. 
The logarithmic Sobolev inequality for the Gaussian measure $\NN_{0,1}$ 
(see e.g. \cite[formula (1.2)]{Gro75}) yields, for 
each $t\in \R$,
\begin{align*}
\int\limits_{\R^n}|u(t,x)|^2\log(|u(t,x)|) \NN_{0,1}(\ddc x)
\leq&
 \int\limits_{\R^n}|\D_x u|^2
\NN_{0,1}(\ddc x)\\
&+\|u(t,\cdot)\|^2_{L^2(\R^n,\NN_{0,1})}\log
\|u(t,\cdot)\|_{L^2(\R^n,\NN_{0,1})}.
\end{align*}
Hence, following e.g. the 
lines of
\cite{LMP06},  for each $k>1$ we obtain
\begin{align*}
\int\limits_0^T\int\limits_{B(0,R)^c}|u|^{2}\NN_{0,1}(\ddc x)\dd t
\leq
&\int\limits_0^T\int\limits_{B(0,R)^c }\chi_E(x) k^2\NN_{0,1}(\ddc 
x)\dd t
\\
&+\frac{1}{\log k}\int\limits_0^T\int\limits_{B(0,R)^c
}\chi_{E^c}(x)|u|^2\log|u|\NN_{0,1}(\ddc x)\dd t
\\
\leq &k^2T \NN_{0,1}(B(0,R)^c )+\frac{T}{\log k},
\end{align*}
where $E=\{|u|<k\}$. Therefore, given $\varepsilon>0$,
there exists $R>0$, independent of $u$, such that
\begin{align*}
\int\limits_0^T\int\limits_{B(0,R)^c}|u|^2\NN_{0,1}(\ddc x)\dd 
t\leq\varepsilon.
\end{align*}
Since $L^2_\#( (0,T)\times B(0,R),\ddc t \times \NN_{0,1} )=L^2_\#(
(0,T)\times B(0,R))$, and the embedding of $H^{1,2}( (0,T)\times
B(0,R))$ into $L^2( (0,T)\times B(0,R) )$ is compact, we find
$\{f_1,\dots,f_k\}\subset L^2( (0,T)\times B(0,R) )$ such that the
balls $B(f_i,\varepsilon)$ cover the restrictions of the functions of
$B$ to $(0,T)\times B(0,R)$. Now, let $\tilde f_i$ denote the
extensions to $(0,T)\times \R^n$ by $0$. Then
$B\subset\cup_{i=1}^kB(\tilde f_i,2\varepsilon)$ and the proof is
complete.
\end{proof}


\section{Proof of Theorem~\ref{thm:maxreg}}
Without loss of generality, we restrict ourselves to the case 
$T_1=a<0$ and $T_2=0$. We
first consider the problem
\begin{align*}
\left\{\begin{array}{l}
u_s(s,x) + \cK(s)u(s,x)=0, \;\;s\in(a,0), \;x\in \R^{n},
\\
\\
u(0,x)=\varphi(x), \;x\in \R^{n}. 
\end{array}\right.
\end{align*}

\begin{proposition}
\label{pr:maxreginit}
For each $\varphi \in H^1(\R^n,\nu_0)$ the function $(s,x)\mapsto 
u(s,x):= (P_{s,0}\varphi)(x)$ belongs to 
$H^{1,2}((a,0)\times\R^n,\nu)$ and there exists $C>0$, independent of 
$\varphi$, such that 
\begin{equation}
\label{stimaPs0}
\|u\|_{ H^{1,2}((a,0)\times\R^n,\nu)} \leq 
C\|\varphi\|_{H^1(\R^n,\nu_0)}.
\end{equation}
\end{proposition}
\begin{proof}
We use the following identities:
\begin{equation}
\label{eq:commutatore}
[\LL(s), 
D]P_{s,0}\varphi=(\LL(s)D-DL(s))P_{s,0}\varphi=A(s)DP_{s,0}\varphi=A(s)U^*(s,0)P_{s,0}D\varphi, 
\end{equation}
 and, for $\psi\in  H^3(\R^n, 
\nu_s)$, 
\begin{equation}
\label{eq:gradiente/derivata}
\int_{\R^n} |D\psi|^2 \partial_s \rho(s,x)\ddc x = 
2 \int_{\R^n} \langle \LL(s) D\psi, D\psi\rangle   \ddc \nu_s  + 
\int_{\R^n} |B^*(s)D^2\psi|^2  \ddc \nu_s . 
\end{equation}
Formula \eqref{eq:commutatore}  follows from the explicit expressions 
of $\LL(s)$ and $P_{s,0}$,
while \eqref{eq:gradiente/derivata} follows from Lemma~\ref{invop} 
and the identity 
$\LL (s)(\varphi^2) = 2 \varphi \LL(s)\varphi + | B^*(s)D\varphi|^2$, 
applied to each derivative $D_{j}\psi$. 

Thus, we obtain 	
\begin{align*}
\partial_s\int\limits_{\R^n}|DP_{s,0}\varphi|^2\nu_s(\ddc	
x)=&-2\int\limits_{\R^n}\langle
D\LL (s)P_{s,0}\varphi,DP_{s,0}\varphi\rangle
\nu_s(\ddc x)
+\int\limits_{\R^n}|DP_{s,0}\varphi|^2
\partial_s \rho(s,x)\dd x\\=&-2\int\limits_{\R^n}\langle 
\LL (s)DP_{s,0}\varphi,
DP_{s,0}\varphi\rangle  \nu_s(\ddc x)\\
&+2\int\limits_{\R^n}\langle [\LL (s),D]P_{s,0}\varphi,
DP_{s,0}\varphi\rangle\nu_s(\ddc x)
	+\int\limits_{\R^n}|DP_{s,0}\varphi|^2
\partial_s\rho(s,x)\dd	x
\\
=&\int\limits_{\R^n}|B^* D^2P_{s,0}\varphi|^2  \nu_s(\ddc x)
\\
&+2\int\limits_{\R^n}\langle
A(s)U^*(s,0)P_{s,0}D\varphi,
U^*(s,0)P_{s,0}D\varphi \rangle
	\nu_s(\ddc x).
\end{align*}
Integrating with respect to $s$, we obtain
\begin{align*}
&\|\,|D\varphi| \,\|^2_{L^2(\R^n,\nu_0)} -
\|\,|DP_{a,0}\varphi|\,\|^2_{L^2(\R^n,\nu_a)}
=
\int\limits_a^0\int\limits_{\R^n}
|DP_{s,0}\varphi|^2 \nu_s(\ddc x)\ddc s
\\
&=\int\limits_{a}^0\int\limits_{\R^n} 
\left( |B^*(s)D^2P_{s,0}\varphi|^2
+2\langle
A(s)U^*(s,0)P_{s,0}D\varphi,
U^*(s,0)P_{s,0}D\varphi
\rangle\right) \nu_s(\ddc x)\ddc s.
\end{align*}
Since $\|B^*(s)^{-1}\| \leq 1/\mu_0$ by assumption \eqref{invertb}, 
then 
$$\|\, |D^2_xu|\,\|_{L^{2}((a,0)\times\R^n, \nu)} \leq 
\frac{1}{\mu_0^2} \int\limits_{a}^0\int\limits_{\R^n} 
|B^*(s)D^2P_{s,0}\varphi|^2
\nu_s(\ddc x)\ddc s$$
and, hence, 
\begin{align*}
\|\, |D^2_xu|\, \|_{L^{2}((a,0)\times\R^n, \nu) }+
\|\,|D_xu(a, \cdot)|\, \|_{L^2(\R^n,\nu_a)}
\leq C(T)
\|\,| D\varphi |\, \|_{H^1(\R^n,\nu_0)},
\end{align*}
where $C(a)>0$ is independent of $\varphi$. Since $\partial_s u = 
-\LL (s)u(s, \cdot)$ the statement follows using estimate 
\eqref{unibound}. \end{proof}

We also need the following lemma about the traces at $t=0$ of 
functions belonging to 
$H^{1,2}((a,0)\times\R^n,\nu)$. 

\begin{lemma}
\label{lem:maxreginit}
We have
$$H^1(\R^n, \nu_0)= \{ u(0, \cdot): \; u\in  
H^{1,2}((a,0)\times\R^n,\nu)\},$$
and the norm
\begin{align*}
 \varphi \mapsto \inf\{ \|  
	u\|_{H^{1,2}( (a,0)\times\R^n, \nu )}:\; u(0, \cdot) = \varphi\}
 \end{align*}
 is equivalent to the norm of $H^1(\R^n, \nu_0)$. 
\end{lemma}
\begin{proof}
By Lemma~\ref{lem:isom}, we have $\TT _0u\in H^{1,2}(
(a,0)\times\R^n, \nu_0) $ and there exists
$C>0$, independent of $u$, such that
\begin{align*}
\|\TT _0
	u\|_{H^{1,2}( (a,0)\times\R^n, \nu_0)}\leq C 
	\| u\|_{H^{1,2}( (a,0)\times\R^n, \nu)}.
\end{align*}
Therefore, by standard arguments,
\begin{align*}
\|u(0,\cdot)\|_{H^1(\R^n)}=\|(\TT _0
u)(0,\cdot)\|_{H^1(\R^n,\nu_0)}
\leq C\|\TT _0 u\|_{H^{1,2}(
(a,0)\times\R^n,\nu_0)}\leq C 
\|u\|_{H^{1,2}(	(a,0)\times\R^n,\nu)},
\end{align*}
where $C$ is independent of $u$.
On the other hand, Proposition \ref{pr:maxreginit} states that for 
each $\varphi \in H^1(\R^n, \nu_0)$ the function $u(s,x) = 
(P_{s,0}\varphi)(x) $ belongs to $H^{1,2}( (a,0)\times\R^n, \nu )$, 
with estimate
\eqref{stimaPs0}. The statement follows. \end{proof}

Finally, we are in the position to prove Theorem~\ref{thm:maxreg}.
Let $f\in L^2( (a,0)\times\R^n,\nu)$, fix $\lambda >0$  and set
$$f_\lambda(s,x)=\left\{ \begin{array}{ll}
-{\rm e}^{\lambda s}f(s,x) & x\in\R^n, \; s\in(a,0), 
\\
\\
0 & x\in\R^n, \; s\notin(a,0)
\end{array}\right. $$
 Then $f_\lambda\in L^2(\R^{1+n},\nu)$ and, by
Theorem~\ref{main1}, $u_\lambda :=(\lambda-G)^{-1}f_\lambda\in
H^{1,2}(\R^{1+n},\nu)$, and 
$\| u_\lambda\|_{H^{1,2}(\R^{1+n},\nu)} \leq C \|f_\lambda\|_{L^{ 
2}(\R^{1+n},\nu)}$ with $C$ independent on $f$. Moreover, 
$u_1(s,x):={\rm e}^{-\lambda
s}u_\lambda(s,x)$, $x\in\R^n$, $s\in(a,0)$ satisfies
\begin{align*}
\partial_s u_1(s,x)+\cL(s)u_1(s,x)=-\lambda{\rm
e}^{-\lambda s}u_\lambda(s,x)+{\rm e}^{-\lambda
s}\cG u_\lambda =-{\rm e}^{-\lambda
s}f_\lambda(s,x)=f(s,x)
\end{align*}
for $x\in\R^n$ and $s\in(a,0)$.
Furthermore, there exists $C_1>0$, independent of $f$,  such that
\begin{align*}
\|u_{1}\|_{H^{1,2}( (a,0)\times\R^n)}\leq C_1\|f\|_{L^2(
(a,0)\times\R^n)}.
\end{align*}
Hence, by Lemma~\ref{lem:maxreginit},
\begin{align*}
u:=u_1+P_{\cdot,0}(u_0-u_1(0,\cdot))\in
H^{1,2}( (a,0)\times\R^n,\nu)
\end{align*}
and $u$ satisfies \eqref{e7} and  \eqref{estmaxreg}.

\end{document}